\documentclass[12pt]{article}
\usepackage[russian]{babel}
\usepackage{amsthm,amsfonts,amssymb,amscd}

\begin{document}

\centerline{\large \bf Birational geometry of Fano direct
products}\vspace{1cm}

\centerline{\bf A.V.Pukhlikov}\vspace{0.3cm}

\noindent {\small Department of Mathematical Sciences, The
University of Liverpool, Liverpool L69 7ZL, England
\newline Steklov Institute of Mathematics, Gubkina 8, 119991
Moscow, Russia
\newline e-mail: pukh@liv.ac.uk, pukh@mi.ras.ru}
\vspace{0.3cm}

\parshape=1
3cm 10cm \noindent {\small \quad\quad\quad \quad\quad\quad\quad
\quad\quad\quad {\bf Abstract}\newline We prove birational
superrigidity of direct products $V=F_1\times\dots\times F_K$ of
primitive Fano varieties of the following two types: either
$F_i\subset{\mathbb P}^M$ is a general hypersurface of degree
$M$, $M\geq 6$, or $F_i\stackrel{\sigma}{\to}{\mathbb P}^M$ is a
general double space of index 1, $M\geq 3$. In particular, each
structure of a rationally connected fiber space on $V$ is given
by a projection onto a direct factor. The proof is based on the
connectedness principle of Shokurov and Koll\' ar and the
technique of hyper\-tangent divisors.}\vspace{1cm}

\section*{Introduction}

\subsection{Fano direct products}
A smooth projective variety $F$ is a {\it primitive Fano
variety}, if  $\mathop{\rm Pic}F={\mathbb Z}K_F$, the
anticanonical class ($-K_F$) is ample and $\mathop{\rm dim}F\geq
3$. Typical examples of primitive Fano varieties are given by
smooth complete intersections  of index 1 in the (weighted)
projective spaces (in dimension 4 there are no other examples).
The most visual examples are Fano hypersurfaces $X_d\subset
\mathbb P^d$, $\mathop{\rm deg} X_d=d$, $d\geq 4$ and Fano double
spaces $\sigma\colon X\to {\mathbb P}^d$, branched over a smooth
hypersurface $W_{2d}\subset {\mathbb P}^d$ of degree $2d$.

It is known (see, for instance, [1,13,14,22-26,28,29,31,33,34]),
that typical primitive Fano varieties (in dimension $\geq 4$
probably all primitive Fano varieties) are {\it birationally
rigid}, in particular, they cannot be fibered by a non-trivial
rational map into rationally connected varieties (for the precise
definition, see below). The aim of this paper is to describe
birational geometry of direct products of the form
$$
V=F_1\times\dots\times F_K,
$$
$K\geq 2$, where $F_i$ are primitive Fano varieties. We will
prove that provided that the direct factors $F_i$  satisfy
certain conditions, all structures of rationally connected
fibrations on $V$ are given by the projections $V\to
F_{i_1}\times\dots\times F_{i_k}$ onto direct factors. More
precisely, the variety $V$ itself is birationally rigid.

Roughly speaking, the variety $V$ can be reconstructed from its
field of rational functions ${\mathbb C}(V)$ (we work over the
ground field of complex numbers ${\mathbb C}$). The conditions
that we have to impose on the direct factors $F_i$ are very
strong. In this paper, we show that sufficiently general  (in the
sense of Zariski topology) Fano hypersurfaces $X_d\subset
{\mathbb P}^d$, $d\geq 6$, and Fano double covers of index 1
satisfy these conditions.

\subsection{Birational rigidity and superrigidity}
For a smooth projective variety $X$ denote by $A^iX$ the group of
classes of cycles of codimension $i$ modulo rational equivalence.
In particular, $A^1X\cong\mathop{\rm Pic}X$ is the Picard group.
Set
$$
A^1_+X\subset A^1X\otimes{\mathbb R}
$$
to be the closed cone of pseudoeffective classes, that is, the
smallest closed cone containing the classes of all effective
divisors. Assume that the variety $X$ is rationally connected
[9,16]. For an effective divisor $D$ on $X$ define its {\it
threshold of canonical adjunction} by the formula
\begin{equation}\label{i1}
c(D)=c(D,X)=\mathop{\rm sup}\{\varepsilon\in{\mathbb
Q}_+|D+\varepsilon K_X \in A^1_+X\},
\end{equation}
where $K_X$ is the canonical class. For an arbitrary non-empty
linear system $\Sigma$ on $X$ we write $c(\Sigma)=c(D)$, where
$D\in \Sigma$ is an arbitrary divisor. For a {\it movable} linear
system $\Sigma$ we define its {\it virtual threshold of canonical
adjunction} by the formula
$$
c_{\rm virt}(\Sigma)=\mathop{\rm inf}\limits_{\widetilde X\to
X}\{c(\widetilde\Sigma,\widetilde X)\},
$$
where the infimum is taken over all sequences of blow ups
$\widetilde X\to X$ with non-singular centres (or, equivalently,
over all birational morphisms of smooth varieties $\widetilde
X\to X$) and $\widetilde \Sigma$ is the strict transform of the
system $\Sigma$ on $\widetilde X$.

{\bf Definition 1.} (i) The variety $X$ is said to be {\it
birationally superrigid}, if for any movable linear system
$\Sigma$ on $X$ the equality
$$
c(\Sigma)=c_{\rm virt}(\Sigma)
$$
holds.

(ii) The variety $X$ is said to be {\it birationally rigid}, if
for any movable linear system $\Sigma$ there exists a birational
self-map $\chi\in\mathop{\rm Bir}X$ such that
$$
c_{\rm virt}(\Sigma)= c(\chi^{-1}_* \Sigma).
$$
Here $\chi^{-1}_* \Sigma$ means the {\it strict transform} of the
system $\Sigma$ with respect to $\chi$.

{\bf Remark 1.} The virtual threshold of canonical adjunction is
birationally invariant in the following sense: let
$$
\chi\colon X_1-\,-\,\to X_2
$$
be a birational map, $\Sigma_i$ moving linear system on $X_i$,
$i=1,2$, and moreover, $\Sigma_1$ the strict transform of
$\Sigma_2$ with respect to $\chi$. Then
$$
c_{\rm virt}(\Sigma_1)=c_{\rm virt}(\Sigma_2),
$$
although the actual thresholds $c(\Sigma_1)$ and $c(\Sigma_2)$
should not be the same.

{\bf Remark 2.} Let $\pi\colon X\to S$ be a rationally connected
fiber space, that is, for a general point $z\in S$ the fiber
$\pi^{-1}(z)$ is a rationally connected variety. Assume that the
base $S$ is non-trivial (that is, $\mathop{\rm dim} S\geq 1$) and
$$
\Sigma=\pi^*\Lambda
$$
is the pull back of a linear system $\Lambda$ on $S$. It is easy
to see that
$$
c(\Sigma)= 0,
$$
since $K_X$ is negative on the fibers of $\pi$. Therefore,
$$
c_{\rm virt}(\Sigma)=c(\Sigma)=0.
$$
Thus if on a rationally connected variety $X$ there is a {\it
structure} of a rationally connected fiber space (that is, a
non-trivial {\it rational} map $\varphi\colon X -\,-\,\to S$, the
general fiber of which is rationally connected), then on $X$
there is a movable linear system $\Sigma$, satisfying the
equality $c_{\rm virt}(\Sigma)=0$. In particular, if the variety
$X$ is birationally rigid and there is a non-trivial structure of
a rationally connected fiber space on $X$, then there is a
movable linear system $\Sigma$ satisfying the equality
$c(\Sigma)=0$. For trivial reasons, on a primitive Fano variety
$F$ there is no such systems. Therefore, birationally rigid
primitive Fano varieties {\it do not admit non-trivial structures
of a rationally connected fiber space.}

The first example of a birationally rigid (in fact, superrigid)
primitive Fano variety was given by the theorem of V.A.Iskovskikh
and Yu.I.Manin on the three-dimensional quartics [14]. The method
of proving it ({\it the method of maximal singularities}) later
accumulated a number of new ideas and was considerably improved.
The word combination ``birational rigidity'' was introduced in
the late 80ies, see Sec. 0.7 below.

\subsection{The main result}

{\bf Definition 2.} We say that a primitive Fano variety $F$ is
{\it divisorially canonical}, or satisfies the condition ($C$)
(respectively, is {\it divisorially log canonical}, or satisfies
the condition ($L$)), if for any effective divisor $D\in|-nK_F|$,
$n\geq 1$, the pair
\begin{equation}\label{i2}
(F,\frac{1}{n}D)
\end{equation}
has canonical (respectively, log canonical) singularities. If the
pair (\ref{i2}) has canonical singularities for a general divisor
$D\in \Sigma \subset|-nK_F|$ of any {\it movable} linear system
$\Sigma$, then we say that $F$ satisfies the condition of {\it
movable canonicity}, or the condition ($M$).

Explicitly, the condition ($C$) is formulated in the following
way: for any birational morphism $\varphi\colon \widetilde F\to
F$ and any exceptional divisor $E\subset \widetilde F$ the
following inequality
\begin{equation}\label{i3}
\nu_E(D)\leq na(E)
\end{equation}
holds. The inequality (\ref{i3}) is opposite to the classical
{\it Noether-Fano} inequality [14,25,26,34]. The condition ($L$)
is weaker: the inequality
\begin{equation}\label{i4}
\nu_E(D)\leq n(a(E)+1)
\end{equation}
is required. In (\ref{i3}) and (\ref{i4}) the number $a(E)$ is
the discrepancy of the exceptional divisor $E\subset \widetilde
F$ with respect to the model $F$. The inequality (\ref{i4}) is
opposite to the {\it log Noether-Fano inequality}. The condition
($M$) means that (\ref{i3}) holds for a general divisor $D$ of
any movable linear system $\Sigma\subset |-nK_F|$ and any
discrete valuation $\nu_E$.

It is well known (essentially starting from the classical paper
of V.A.Iskovskikh and Yu.I.Manin [14]) that the condition ($M$)
ensures birational superrigidity. This condition is proved for
many classes of primitive Fano varieties, see [14,22,23,28,31,34].
Note also that the condition ($C$) is stronger than both ($L$) and
($M$).

Now let us formulate the main result of the present paper.

{\bf Theorem 1.} {\it Assume that primitive Fano varieties
$F_1,\dots,F_K$, $K\geq 2$, satisfy the conditions ($L$) and
($M$). Then their direct product
$$
V=F_1\times\dots\times F_K
$$
is birationally superrigid.}

{\bf Remark 3.} For $K=1$ the claim of the theorem is obviously
true (the condition ($M$) suffices).

{\bf Corollary 1.} {\it {\rm (i)} Every structure of a rationally
connected fiber space on the variety $V$ is given by a projection
onto a direct factor. More precisely, let $\beta\colon
V^{\sharp}\to S^{\sharp}$ be a rationally connected fiber space
and $\chi\colon V-\,-\,\to V^{\sharp}$ a birational map. Then
there exists a subset of indices
$$
I=\{i_1,\dots,i_k\}\subset \{1,\dots,K\}
$$
and a birational map
$$
\alpha\colon F_I=\prod\limits_{i\in I}F_i-\,-\,\to S^{\sharp},
$$
such that the diagram
$$
\begin{array}{rcccl}
& V &\stackrel{\chi}{-\,-\,\to} & V^{\sharp}&\\
\pi_I &\downarrow & &\downarrow &\beta\\
& F_I & \stackrel{\alpha}{-\,-\,\to}& S^{\sharp}&
\end{array}
$$
commutes, that is, $\beta\circ\chi=\alpha\circ\pi_I$, where
$\pi_I\colon\prod\limits^K_{i=1}F_i\to \prod\limits_{i\in I}F_i$
is the natural projection onto a direct factor. In particular,
the variety $V$ admits no structures of a fibration into
rationally connected varieties of dimension smaller than
$\mathop{\rm min}\{\mathop{\rm dim}F_i\}$. In particular, $V$
admits no structures of a conic bundle or a fibration into
rational surfaces.

{\rm (ii)} The groups of birational and biregular self-maps of
the variety $V$ coincide:
$$
\mathop{\rm Bir}V=\mathop{\rm Aut}V.
$$
In particular, the group $\mathop{\rm Bir}V$ is finite.

{\rm (iii)} The variety $V$ is non-rational.}

{\bf Remark 4.} The group of biregular automorphisms $\mathop{\rm
Aut}V$ is easy to compute. Let us break the set $F_1,\dots,F_K$
into subsets of pair-wise isomorphic varieties:
$$
I=\{1,\dots,K\}= \mathop{\bigvee}\limits^l_{k=1}I_k,
$$
where $F_i\cong F_j$ if and only if $\{i,j\}\subset I_k$ for some
$k\in\{1,\dots,l\}$. It is easy to see that
$$
\mathop{\rm Aut}V=\prod\limits^l_{j=1}\mathop{\rm
Aut}(\prod\limits_{i\in I_j}F_i).
$$
In particular, if the varieties $F_1,\dots,F_K$ are pair-wise
non-isomorphic, we get
$$
\mathop{\rm Aut}V=\prod\limits^K_{i=1}\mathop{\rm Aut}F_i
$$
(and this group acts on $V$ component-wise). In the opposite
case, if
$$
F_1\cong F_2\cong \dots\cong F_K\cong F,
$$
we obtain the exact sequence
$$
1\to(\mathop{\rm Aut}F)^{\times K}\to \mathop{\rm Aut}V \to S_K
\to 1,
$$
where $S_K$ is the symmetric group of permutations of $K$
elements. In fact, in this case $\mathop{\rm Aut}V$ contains a
subgroup isomorphic to $S_K$ which permutes direct factors of
$V$, so that $\mathop{\rm Aut}V$ is a semi-direct product of the
normal subgroup $(\mathop{\rm Aut}F)^{\times K}$ and the
symmetric group $S_K$.

It seems that the following generalization of Theorem 1 is true.

{\bf Conjecture 1.} {\it Assume that $F_1,\dots,F_K$ are
birationally (super)rigid primitive Fano varieties. Then their
direct product $V=F_1\times\dots\times F_K$ is birationally
(super)rigid.}

\subsection{Divisorially canonical varieties}

Let ${\mathbb P}={\mathbb P}^M$, $M\geq 3$, be the complex
projective space. Set
$$
{\cal F}={\mathbb P}(H^0({\mathbb P},{\cal O}_{\mathbb
P}(M))),\quad {\cal W}={\mathbb P}(H^0({\mathbb P},{\cal
O}_{\mathbb P}(2M)))
$$
be the spaces of hypersurfaces of degrees $M$ and $2M$,
respectively.

{\bf Theorem 2.} {\it {\rm (i)} For $M\geq 6$ there exists a
non-empty Zariski open subset ${\cal F}_{\rm reg}\subset {\cal
F}$ such that any hypersurface $F\in {\cal F}_{\rm reg}$ is
non-singular and satisfies the condition ($C$).

{\rm (ii)} There exists a non-empty Zariski open subset ${\cal
W}_{\rm reg}\subset {\cal W}$ such that for any $W\in {\cal
W}_{\rm reg}$ the Fano double space $\sigma\colon F\to {\mathbb
P}$, branched over the hypersurface $W$, satisfies the condition
($C$).}

The sets of regular hypersurfaces ${\cal F}_{\rm reg}$, ${\cal
W}_{\rm reg}$ are explicitly described in Sec. 2.1 and 2.2, and
their non-emptiness is proved in Sec. 2.3. It seems that the
following generalization of Theorem 2 is true.

{\bf Conjecture 2.} {\it A general element $F\in \Phi$ of any
family of primitive Fano complete intersec\-tions $\Phi$ in a
weighted projective space satisfies the condition ($C$). In
particular, this condition is satisfied for a general quartic
$F=F_4\subset {\mathbb P}^4$ and a general quintic $F=F_5\subset
{\mathbb P}^5$.}

\subsection{The scheme of the proof and the structure of the paper}

Theorem 1 is proved by induction on the number of factors $K$.
For $K=1$, as we have mentioned above, the claim of the theorem
holds in a trivial way. So assume that $K\geq 2$.

Assume the converse: there is a moving linear system $\Sigma$ on
$V$ such that the inequality
$$
c_{\rm virt}(\Sigma) < c(\Sigma)
$$
holds. By the definition of the virtual threshold of canonical
adjunction it means that there exists a sequence of blow ups
$\varphi\colon \widetilde V\to V$ such that the inequality
\begin{equation}\label{i5}
c(\widetilde\Sigma)< c(\Sigma)
\end{equation}
holds, where $\widetilde\Sigma$ is the strict transform of the
linear system $\Sigma$ on $\widetilde V$.

To prove Theorem 1, we must show that the inequality (\ref{i5})
is impossible, that is, to obtain a contradiction.

Let $H_i=-K_{F_i}$ be the positive generator of the group
$\mathop{\rm Pic}F_i$. Set
$$
S_i=\prod\limits_{j\neq i}F_i,
$$
so that $V\cong F_i\times S_i$. Let $\rho_i\colon V\to F_i$  and
$\pi_i\colon V\to S_i$ be the projections onto the factors.
Abusing our notations, we write $H_i$ instead of $\rho^*_i H_i$,
so that
$$
\mathop{\rm Pic}V=\mathop{\bigoplus}\limits^K_{i=1}{\mathbb Z}H_i
$$
and
$$
K_V=-H_1-\dots-H_K.
$$
We get
$$
\Sigma\subset|n_1H_1+\dots+n_KH_K|,
$$
whereas
$$
c(\Sigma)=\mathop{\rm min}\{n_1,\dots,n_K\}.
$$
Without loss of generality we assume that $c(\Sigma)=n_1$. By the
inequality (\ref{i5}) we get $n_1\geq 1$. Set $n=n_1$,
$\pi=\pi_1$, $F=F_1$, $S=S_1$. We get
$$
\Sigma \subset |-nK_V+\pi^*Y|,
$$
where $Y=\sum\limits^K_{i=2}(n_i-n)H_i$ is an effective class on
the base $S$ of the fiber space $\pi$.

Now we need to modify the birational morphism $\varphi$. Here our
arguments are similar to the proof of the Sarkisov theorem
[38,39]. For an arbitrary sequence of blow ups $\mu_S\colon
S^+\to S$ we set $V^+=F\times S^+$ and obtain the following
commutative diagram:
\begin{equation}\label{i6}
\begin{array}{rcccl}
&V^+&\stackrel{\mu}{\to}&V&\\
\pi_+&\downarrow&&\downarrow&\pi\\
&S^+&\stackrel{\mu_S}{\to}&S,
\end{array}
\end{equation}
where $\pi_+$ is the projection and $\mu=(\mathop{\rm
id}_F,\mu_S)$. Let $E_1,\dots,E_N\subset \widetilde V$ be all
exceptional divisors of the morphism $\varphi$.

{\bf Proposition 1.} {\it There exists a sequence of blow ups
$\mu_S\colon S^+\to S$ such that in the notations of the diagram
(\ref{i6}) the centre of each discrete valuation $E_i$,
$i=1,\dots,N$, covers either $S^+$ or a divisor on $S^+$:}
$$
\mathop{\rm codim} [\pi_+(\mathop{\rm centre}(E_i,V^+))]\leq 1.
$$

Let $\Sigma^+$ be the strict transform of the linear system
$\Sigma$ on $V^+$. Now the arguments break into two parts due to
the following fact.

{\bf Proposition 2.} {\it The following alternative holds:

{\rm (i)} either the inequality $c(\Sigma^+)< c(\Sigma)$ is true,

{\rm (ii)} or for a general divisor $D^+\in \Sigma^+$ the pair
$$
(V^+,\frac{1}{n}D^+)
$$
is not canonical, and moreover, for some $i=1,\dots,N$ the
discrete valuation $E_i$ determines a non-canonical singularity
of this pair.}

{\bf Remark 5.} The alternative of Proposition 2 should be
understood in the ``and/or'' sense: {\it at least one} of the two
possibilities (i), (ii) takes place (or both).

Now the proof of Theorem 1 is completed in the following way. If
the situation (i) takes place, then we show that the variety $S$
is not birationally superrigid, which contradicts to the
induction hypothesis, since $S$ is a direct product of $K-1$ Fano
varieties. On the other hand, if the case (ii) holds, then by
inversion of adjunction, which follows from the Shokurov-Koll\'
ar connectedness principle, we obtain a contradiction with the
condition ($L$) for the fiber $F$. Now the proof of Theorem 1 is
complete.

The scheme of arguments described above is realized in full
detail in Sec. 1. The subsequent structure of the paper is as
follows. In Sec. 2 we prove Theorem 2. In Sec. 2.1 we describe
the regularity conditions, defining the open set ${\cal F}_{\rm
reg}$ and prove divisorial canonicity of hypersurfaces
$F_M\subset {\mathbb P}^M$, $M\geq 6$, $F_M\in {\cal F}_{\rm
reg}$. In Sec. 2.2 we describe the regularity conditions,
defining the open set  ${\cal W}_{\rm reg}$ and prove divisorial
canonicity of double spaces $F\to{\mathbb P}^M$, branched over
hypersurfaces $W\in {\cal W}_{\rm reg}$. In Sec. 2.3 it is proved
that the sets ${\cal F}_{\rm reg}$, ${\cal W}_{\rm reg}$ are
non-empty (for ${\cal F}_{\rm reg}$ it is not quite obvious). For
the hypersurfaces the proof is based on the {\it technique of
hypertangent divisors} [28,31,34].

The scheme of our proof of Theorem 2 is as follows. Assume that
the pair $(F,\frac{1}{n}D)$ is not canonical, where $D\in
|-nK_F|$, $n\geq 1$, $F$ is a regular Fano hypersurface or a
double space. It is easy to show that the centre of a
non-canonical singularity is a point $x\in F$. By inversion of
adjunction one obtains

{\bf Proposition 3.} {\it Let $\varphi\colon \widetilde F\to F$
be the blow up of the point $x$, $E\subset \widetilde F$ the
exceptional divisor. For some hyperplane $B\subset E$ the
inequality
\begin{equation}\label{i7}
\mathop{\rm mult}\nolimits_xD+\mathop{\rm
mult}\nolimits_B\widetilde D> 2n
\end{equation}
holds, where $\widetilde D\subset \widetilde F$ is the strict
transform of the divisor $D$.}

If $F$ is a double space, then it is not hard to check that the
inequality (\ref{i7}) is impossible for an effective divisor
$D\in|-nK_F|$. If $F\subset{\mathbb P}^M$ is a Fano hypersurface,
the contradiction is obtained by means of a (modified) technique
of hypertangent divisors [28,31,34].

Section 3 contains a proof of Proposition 3. There we also
discuss briefly the connectedness principle and inversion of
adjunction in the form in which we need these facts. For a
detailed exposition of these fundamental results, see the
original work [15,40]. We include a brief explanation of these
facts in this paper because of their importance.

\subsection{Remarks}

(i) Repeating the proof of Theorem 1 word for word, one can
obtain a stronger fact, that is, the coincidence $c_ {\rm
virt}(D)=c(D)$ of the virtual and actual thresholds of canonical
adjunction for an arbitrary {\it irreducible} divisor $D$ on $V$.
The irreducibility is essential. It is in this way that we argued
in the first version of this paper [37]. However, in [37] instead
of the termination of canonical adjunction in the sense of the
formula (\ref{i1}) we used a formally stronger condition that for
some $\varepsilon\in{\mathbb Q}$ and a sequence of blow ups
$\widetilde V\to V$ the class $\widetilde D+\varepsilon
K_{\widetilde V}$ is negative on the curves from some family
sweeping out $\widetilde V$. This approach does not give
birational superrigidity, but makes it possible to describe
(exactly in the same way as in this paper) all structures of a
rationally connected fiber space on $V$. In particular, it gives
all the claims of Corollary 1. Proof of Theorem 2 is completely
the same as in [37].

(ii) Let us emphasize some points where the present paper differs
from the previous papers on this subject. It is natural to call
the method we use here {\it linear}, in contrast to the {\it
quadratic} method of the previous papers: we study singularities
of a {\it single} divisor $D$, whereas in the above-mentioned
papers the main step was to construct the {\it self-intersection}
of a moving linear system $\Sigma$, that is, the effective cycle
$Z=(D_1\circ D_2)$ of codimension two, where $D_i\in \Sigma$ are
general divisors.

Another point where the present paper differs from the previous
ones is that here we for the first time study Fano fiber spaces
with the base and the fiber of arbitrary dimension: in Sarkisov's
papers [38,39] the fiber is one-dimensional, in
[10-12,27,30,35,36,41,42] the base is one-dimensional, that is,
${\mathbb P}^1$.

Finally, note that the starting point of the previous papers was
a condition of ``sufficient twistedness'' of a Fano fiber space
over the base (the Sarkisov condition for conic bundles, the
$K^2$-condition for Fano fiber spaces over ${\mathbb P}^1$),
whereas direct products are not twisted over the base at all.

(iii) The technique developed in this paper opens new ways of
studying movable linear systems on Fano fiber spaces $V/{\mathbb
P}^1$. If $F=F_t$, $t\in {\mathbb P}^1$, the fiber over the point
$t$, satisfies the condition ($L$), then the linear system
$$
\Sigma\subset |-nK_V+lF|
$$
with $l\in {\mathbb Z}_+$ does not admit a maximal singularity
$E$, the centre of which $B=\mathop{\rm centre}(E,V)$ is
contained in the fiber $F$. This remark demonstrates the power of
the linear method compared to the quadratic one. On the other
hand, it is not easy to prove the condition ($L$). The method,
developed in this paper for that purpose, requires stronger
conditions of general position than those used for the quadratic
method.

Another difference of the linear method from the quadratic one is
that it applies to the trivial fibrations $V=F\times{\mathbb
P}^1$, where $F$ is a primitive Fano variety. Further improvement
of the ideas developed in this paper makes it possible to prove
that if $F$ satisfies the condition ($C$), then on $V$ there is
only one conic bundle structure, given by the projection $V\to
F$, whereas all other non-trivial structures of a rationally
connected fiber space are $\mathop{\rm Bir}V$-equivalent to the
fibration $F\times{\mathbb P}^1\to {\mathbb P}^1$ (that is, to
the projection onto the direct factor). The group $\mathop{\rm
Bir}V$ up to the finite group $\mathop{\rm Aut}F\subset
\mathop{\rm Aut}V$ coincides with the group $\mathop{\rm
Aut}{\mathbb P}^1_{{\mathbb C}(F)}=PGL(2,{\mathbb C}(F))$ of the
fractional-linear transformations over the field of rational
functions ${\mathbb C}(F)$. The proof will be given elsewhere.

\subsection{Historical remarks and acknowledgements}

The phenomenon of birational rigidity was guessed by Fano [5-7],
when he tried to study birational geometry of three-folds by
analogy with surfaces. Fano used termination of canonical
adjunction which already showed itself as a crucial point in the
proof of the Noether theorem on the Cremona group $\mathop{\rm
Bir}{\mathbb P}^2$ and in the theorems of Italian school (for
instance, the Castelnuovo rationality criterion). A preimage of
the concept of maximal singularity is also present in Fano's
papers. However, Fano failed to obtain any more or less complete
proofs.

In the modern age of algebraic geometry the outlines of the
theory of birational rigidity were drawn in the papers of
Yu.I.Manin on surfaces over non-closed fields, see [18-20], where,
in particular, it was proved (using earlier results of B.Segre)
that del Pezzo surfaces of degree 1 with the Picard group
${\mathbb Z}$ are birationally superrigid (in the modern
terminology) and cubic surfaces with the Picard group ${\mathbb
Z}$ are birationally rigid, and the groups of birational
self-maps of these surfaces were computed. In the papers [19,20]
an oriented graph was associated with a finite sequence of blow
ups. Its combinatorial invariants are of extreme importance for
the technique of counting multiplicities [26,28,32,34]. The
construction of this graph appears in almost all papers on this
subject. It is used in the present paper, too, when we study a
non log canonical singularity (Sec. 3.3).

The breakthrough into three-dimensional birational geometry was
made in the paper of V.A.Iskovskikh and Yu.I.Manin [14], where
they proved (in the modern terminology) that three-dimensional
quartics in ${\mathbb P}^4$ are birationally superrigid. Formally
speaking, it is proved in [14] that a birational map of a smooth
quartic onto another one is a biregular isomorphism, however the
arguments of that paper do need any modification for proving
birational superrigidity. In [14] it is pointed out that the
method of the paper (without any improvements) gives a similar
result for smooth double spaces of index 1 and dimension three. A
slightly more sophisticated technique (``untwisting'' of maximal
lines) gives birational rigidity of smooth double quadrics of
index 1 and dimension three (however, they are not superrigid,
the group of birational self-maps is fairly big). Detailed proofs
for these two classes of Fano three-folds were published in [13].

The next step was made by Sarkisov in [38,39]. Using the results
of Iskovskikh on the surfaces with a pencil of rational curves,
Sarkisov proved that the structure of a conic bundle is unique if
the discriminant divisor is sufficiently big (``the Sarkisov
condition''). Again, in fact it was birational rigidity that he
proved, but the very concept did not exist yet.

The first attempt to extend the three-dimensional technique of
Iskovskikh and Manin into arbitrary dimension was made by the
author of this paper in [22,23]. In the same years, birational
geometry of a three-dimensional quartic with a non-degenerate
double point was described [24]. However, the proofs were turning
more and more complicated and difficult to understand. The
methods needed to be improved. Besides, by mid-80ies it was
already clear, that for all those varieties that were
successfully studied (quartics, quintics, double spaces and
double quadrics, conic bundles) one and the same property was
proved. However, that property was not explicitly formulated. The
question was, what sort of property it was?

Trying to give an answer to this question, the author of the
present paper introduced the concepts ``birational rigidity'' and
``birational superrigidity''. The first, rather awkward,
definitions of these concepts were given in the author's talk at
a conference at Warwick university in 1991 and published in [25].
These definitions were subsequently modified several times; quite
probably, the current version (Definition 1 above) is not final,
either. In the end of 90ies Corti and Reid introduced another
version of these concepts, based on the geometric properties of
birationally (super)rigid varieties [1-3].

During the last decade, the techniques of proving birational
(super)rigidity have been radically improved. First of all, the
method of hypertangent divisors made it possible to work with
Fano varieties of arbitrarily high degree [28,29,31,34].
Furthermore, in [1] Corti suggested to use inversion of
adjunction, which is based on the Shokurov-Koll\' ar
connectedness principle (which, in its turn, is based on the
Kawamata-Viehweg vanishing theorem), for studying
self-intersection of movable linear systems, as an alternative to
the technique of counting multiplicities (developed in [26] on the
basis of the test class method [14]). Thus in [1] inversion of
adjunction was first used in the theory of birational rigidity.
This approach turned out to be quite effective for certain types
of maximal singularities of linear systems. Finally, the simple
idea of breaking the self-intersection of a linear system into
the horizontal and vertical parts [27] opened the way to a
systematic study of Fano fiber spaces over ${\mathbb P}^1$
[10-12,30,35,36,41,42].

As a result, the method of studying birational geometry of
rationally connected varieties via investigating singularities of
the corresponding movable linear systems, which we for the
traditional reasons call the method of maximal singularities,
nowadays works for a large class of higher-dimensional Fano
varieties and Fano fiber spaces and hopefully will make the basis
of birational classification of higher-dimensional rationally
connected varieties.

The main result of this paper was obtained by the author in
October 2003 during a research stay at Max-Planck-Institut f\" ur
Mathematik in Bonn. The first version of the paper, [37], was
completed in April 2004 in Liverpool. The author is very grateful
to the Max-Planck-Institut for the excellent conditions of work
and hospitality.


\section{Birationally rigid direct products}

In this section we prove Theorem 1.

\subsection{Resolution of singularities}
Let us prove Proposition 1. Let $E\subset \widetilde V$ be the
exceptional divisor of the birational morphism $\varphi\colon
\widetilde V\to V$, $B=\varphi(E)$ the centre of the discrete
valuation $E$ on $V$. Assume that
$$
\mathop{\rm codim}\nolimits_S\pi(B)\geq 2.
$$
Construct a sequence of commutative diagrams
\begin{equation}\label{a1}
\begin{array}{rcccl}
&V_j&\stackrel{\varepsilon_j}{\to}&V_{j-1}&\\
\pi_j&\downarrow&&\downarrow&\pi_{j-1}\\
&S_j&\stackrel{\lambda_j}{\to}&S_{j-1}&,
\end{array}
\end{equation}
where $j=1,\dots,l$, satisfying the following conditions:

(1) $V_0=V$, $S_0=S$, $\pi_0=\pi$;

(2) $V_j=F\times S_j$, $\pi_j$ is the projection onto the factor
$S_j$, $\varepsilon_j=(\mathop{\rm id}_F,\lambda_j)$ for all
$j\geq 1$;

(3) $\lambda_j$ is the blow up of the irreducible subvariety
$$
B_{j-1}=\pi_{j-1}(\mathop{\rm centre}(E,V_{j-1}))\subset S_{j-1},
$$
where $\mathop{\rm codim} B_{j-1}\geq 2$.

It is obvious that the properties (1)-(3) determine the sequence
(\ref{a1}).

{\bf Lemma 1.} {\it The following inequality holds:} $l\leq
a(E,V)$.

{\bf Proof.} Let $\Delta_j\subset V_j$ be the exceptional divisor
of the morphism $\varepsilon_j$. By construction we get
$$
\mathop{\rm centre}(E,V_j)\subset \Delta_j,
$$
so that
$$
\nu_E(\Delta_j)\geq 1.
$$
Now we obtain
$$
a(E,V)=a(E,V_l)+\sum\limits^l_{j=1}\nu_E(\Delta_j)a(\Delta_j,V)\geq
l.
$$
Q.E.D. for the lemma.

Therefore the sequence of diagrams (\ref{a1}) is finite: we may
assume that $\mathop{\rm centre}(E,V_l)$ covers a divisor on the
base $S_l$.

From this fact (by the Hironaka theorem on the resolution of
singularities) Proposition 1 follows immediately.

{\bf Proof of Proposition 2.} We use the notations of Sec. 0.5.
Consider the diagram (\ref{i6}).

Let
$$
\tau\colon V^{\sharp}\to V^+
$$
be the resolution of singularities of the composite map
$$
V^+\stackrel{\mu}{\to}V\stackrel{\varphi^{-1}}{\to}\widetilde V.
$$
Set $\psi=\varphi^{-1}\circ\mu\circ\tau\colon V^{\sharp}\to
\widetilde V$. There exist an open set $U\subset V^{\sharp}$ and
a closed set of codimension two $Y\subset \widetilde V$ such that
$$
\psi_U=\psi|_U\colon U\to \widetilde V\setminus Y
$$
is an isomorphism. Obviously, if $E\subset V^{\sharp}$ is an
exceptional divisor of the morphism $\tau$ and $E\cap
U\neq\emptyset$, then
$$
E\cap U=\psi^{-1}_U(E_i)
$$
for some exceptional divisor $E_i$ of the morphism $\varphi$.

Let $\Sigma^+$ and $\Sigma^{\sharp}$ be the strict transforms of
the linear system $\Sigma$ on $V^+$ and $V^{\sharp}$,
respectively, $\Sigma_U=\Sigma^{\sharp}|_U$. If $D^{\sharp}\in
\Sigma^{\sharp}$ is a general divisor, then
$$
\widetilde D=\psi_U(D^{\sharp}_U)\in\widetilde\Sigma
$$
is a general divisor of the linear system $\widetilde\Sigma$ (we
make no difference between $\widetilde\Sigma$ and its restriction
onto $\widetilde V\setminus Y$, since the set $Y$ is of
codimension two). We know that
$$
\widetilde D+ nK_{\widetilde V}\notin A^1_+\widetilde V,
$$
see (\ref{i5}). Therefore,
\begin{equation}\label{a2}
D^{\sharp}_U+nK_U\notin A^1_+U.
\end{equation}
Let ${\cal E}$ be the set of exceptional divisors of the morphism
$\tau$ with a non-empty intersection with $U$. By (\ref{a2}) we
get
$$
\tau^*(D^++nK^+)|_U-\sum\limits_{E\in{\cal
E}}(\nu_E(D^+)-na^+(E))E_U\notin A^1_+U,
$$
where $K^+$ is the canonical class of $V^+$, $a^+(E)=a(E,V^+)$.
Consequently, either
$$
D^++nK^+\notin A^1_+V^+,
$$
and we are in the case (i) of Proposition 2, or there exists an
exceptional divisor $E\in{\cal E}$, satisfying the Noether-Fano
inequality
$$
\nu_E(D^+)>n\cdot a^+(E),
$$
that is, the discrete valuation $E$ realizes a non-canonical
singularity of the pair ($V^+,\frac{1}{n}D^+)$. In the latter
case we get part (ii) of the alternative of Proposition 2, since
$E\in{\cal E}$ and thus $E=E_i$ for some $i=1,\dots,N$ (as
discrete valuations). Q.E.D. for Proposition 2.

\subsection{Reduction to the base}

Assume that the case ($i$) of the alternative of Proposition 2
takes place, that is, $D^++nK^+\notin A^1_+V^+$. Let $z\in F$ be
a point of general position. Set
$$
S^+_z=\{z\}\times S^+, \,\,S_z=\{z\}\times S.
$$
It is clear that $K^+_z=K^+|_{S^+_z}$ and $K_z=K_V|_{S_z}$ are
the canonical classes $K^+_S=K_{S^+}$ and $K_S$, respectively. Let
$$
\Sigma_z=\Sigma|_{S_z}  \quad\mbox{and}\quad
\Sigma^+_z=\Sigma^+|_{S^+_z}
$$
be the restriction of the linear systems $\Sigma$, $\Sigma^+$
onto $S_z$ and $S^+_z$. Take general divisors $D_z\in\Sigma_z$ and
$D^+_z\in\Sigma^+_z$. We get a movable linear system $\Sigma_z$
on the variety $S=F_2\times\dots\times F_K$. Moreover,
$$
\Sigma_z\subset |n_2H_2+\dots+n_KH_K|,
$$
so that $c(\Sigma_z)=\mathop{\rm min}\{n_2,\dots,n_K\}\geq
n=c(\Sigma$).

{\bf Lemma 2.} {\it The following estimate holds:}
$$
D^++nK^+=\pi^*_+(D^+_z+nK^+_z)
$$

{\bf Proof.} Set ${\cal E}_S$ to be the set of exceptional
divisors of the morphism $\mu_S$. The exceptional divisors of the
morphism $\mu$ are $F\times E=\pi^*_+E$ for $E\in{\cal E}_S$. We
get
$$
K^+_S=\mu^*_SK_S+\sum\limits_{E\in{\cal E}_S}a_EE
$$
and
$$
K^+=\mu^*K_V+\pi^*_+(\sum\limits_{E\in {\cal E}_S}a_E E),
$$
where $a_E=a(E)$ is the discrepancy of the divisor $E$. For some
numbers $b_E\geq 0$ we get
$$
D^+=\mu^*D-\sum\limits_{E\in{\cal E}_S}b_E\pi^*_+E,
$$
whereas for a point $z\in F$ of general position
$$
D^+_z=\mu^*_SD_z-\sum\limits_{E\in{\cal E}_S}b_EE.
$$
Now taking into account that $D+nK_V=\pi^*Y$ and
$D_z+nK_z=D_z+nK_S=Y$, we obtain the claim of the lemma.

{\bf Corollary 2.} $D^+_z+nK^+_z\notin A^1_+S^+$

{\bf Proof.} Indeed, it is clear that
$$
\pi^*_+A^1_+S^+\subset A^1_+V^+.
$$
Q.E.D. for the corollary.

Thus for the strict transform $\Sigma^+_z$ of the linear system
$\Sigma_z$ on $S^+$ we get the inequality
$$
c(\Sigma^+_z)<c(\Sigma_z).
$$
The more so, $c_{\rm virt}(\Sigma_z)<c(\Sigma_z)$. Therefore the
variety $S$ is not birationally superrigid. This contradicts the
induction hypothesis.

\subsection{Reduction to the fiber}

By Proposition 2 and Sec. 1.2  for a general divisor
$D^+\in\Sigma^+$ the pair $(V^+,\frac{1}{n}D^+)$ is not
canonical, that is, there exists a birational morphism
$V^{\sharp}\to V^+$ and an exceptional divisor $E\subset
V^{\sharp}$, satisfying the Noether-Fano inequality
$\nu_E(\Sigma^+)> n\cdot a^+_E$, where $a^+_E=a(E,V^+)$.
Moreover, we can assume that the centre $B=\mathop{\rm
centre}(E,V)$ of the valuation $E$ covers a divisor on the base
or the whole base: $\mathop{\rm codim}_{S^+}T\leq 1$, where
$T=\pi_+(B)$.

Let $t\in T$ be a point of general position. The fiber
$F_t=\pi^{-1}_+(t)$ cannot lie entirely in the base set
$\mathop{\rm Bs}\Sigma^+$ of the moving linear system $\Sigma^+$,
since
$$
\mathop{\rm codim}\nolimits_{V^+}\pi^{-1}_+(T)\leq 1.
$$
Therefore, $\Sigma^+_t=\Sigma^+|_{F_t}$ is a non-empty linear
system on $F$, $\Sigma^+_t\subset |nH|=|-nK_F|$ (if $T\subset
S^+$ is a divisor, then $\Sigma^+_t$ can have fixed components).
Let $D^+_t\in\Sigma^+_t$ be a general divisor. By inversion of
adjunction, see [15] and Sec. 3.2 below, the pair
$$
(F,\frac{1}{n}D^+_t)
$$
is not log canonical. We get a contradiction with the condition
($L$). This contradiction completes the proof of Theorem 1.

\subsection{The structures of a rationally connected fiber space}

Let us prove the claim (i) of Corollary 1. Let $\beta\colon
V^{\sharp}\to S^{\sharp}$ be a rationally connected fiber space,
$\chi\colon V-\,-\,\to V^{\sharp}$ a birational map. Take a very
ample linear system $\Sigma^{\sharp}_S$ on the base $S^{\sharp}$
and let
$$
\Sigma^{\sharp}=\beta^*\Sigma^{\sharp}_S
$$
be a movable linear system on $V^{\sharp}$. As we have mentioned
above (Remark 2), $c(\Sigma^{\sharp})=0$. Let $\Sigma$ be the
strict transform of the system $\Sigma^{\sharp}$ on $V$. By our
remark,
$$
c_{\rm virt}(\Sigma)=0,
$$
so that by Theorem 1 we conclude that $c(\Sigma)=0$. Therefore,
in the presentation
$$
\Sigma\subset |-n_1H_1-\dots-n_KH_K|
$$
we can find a coefficient $n_e=0$. We may assume that $e=1$.
Setting $S=F_2\times\dots\times F_K$ and $\pi\colon V\to S$ to be
the projection, we get
$$
\Sigma\subset |\pi^*Y|
$$
for a non-negative class $Y$ on $S$. But this means that the
birational map $\chi$ of the fiber space $V/S$ onto the fiber
space $V^{\sharp}/S^{\sharp}$ is fiber-wise: there exists a
rational dominant map
$$
\gamma\colon S-\,-\,\to S^{\sharp},
$$
making the diagram
$$
\begin{array}{rcccl}
&V &\stackrel{\chi}{-\,-\,\to} &V^{\sharp}&\\
\pi&\downarrow&&\downarrow&\beta\\
&S&\stackrel{\gamma}{-\,-\,\to}&S^{\sharp}
\end{array}
$$
commutative. For a point $z\in S^{\sharp}$ of general position
let $F^{\sharp}_z=\beta^{-1}(z)$ be the corresponding fiber,
$F^{\chi}_z\subset V$ its strict transform with respect to
$\chi$. By assumption, the variety  $F^{\chi}_z$ is rationally
connected. On the other hand,
$$
F^{\chi}_z=\pi^{-1}(\gamma^{-1}(z))=F\times \gamma^{-1}(z),
$$
where $F=F_1$ is the fiber of $\pi$. Therefore, the fiber
$\gamma^{-1}(z)$ is also rationally connected.

Thus we have reduced the problem of description of rationally
connected structures on $V$ to the same problem for $S$. Now the
claim (i) of Corollary 1 is easy to obtain by induction on the
number of direct factors $K$. For $K=1$ it is obvious that there
are no non-trivial rationally connected structures (Remark 2).
The second part of the claim (i) (about the structures of conic
bundles and fibrations into rational surfaces) is obvious since
$\mathop{\rm dim} F_i\geq 3$ for all $i=1,\dots,K$.
Non-rationality of $V$ is now obvious.

\subsection{Birational self-maps}
Let us prove the claim (ii) of Corollary 1. Set ${\cal RC}(V)$ to
be the set of all structures of a rationally connected fiber
space on $V$ with a non-trivial base. By the part (i) we have
$$
{\cal RC}(V)=\{\pi_I\colon V\to F_I=\prod_{i\in
I}F_i\,|\,\emptyset \neq I\subset \{1,\dots,K\}\}.
$$
The set ${\cal RC}(V)$ has a natural structure of an ordered set:
$\alpha\leq \beta$ if $\beta$ factors through $\alpha$. Obviously,
$$
\pi_I\leq \pi_J \quad \mbox{if and only if}\quad J\subset I.
$$
For $I=\{1\dots,K\}\setminus\{e\}$ set $\pi_I=\pi_e$, $F_I=S_e$.
It is obvious that $\pi_1\dots,\pi_K$ are the minimal elements of
${\cal RC}(V)$.

Let $\chi\in \mathop{\rm Bir}V$ be a birational self-map. The map
$$
\chi^*\colon{\cal RC}(V)\to {\cal RC}(V),
$$
$$\chi^*\colon \alpha\longmapsto \alpha \circ\chi,
$$
is a bijection preserving the relation $\leq$. From here it is
easy to conclude that $\chi^*$ is of the form
$$
\chi^*\colon \pi_I\longmapsto \pi_{I^{\sigma}},
$$
where $\sigma\in S_K$ is a permutation of $K$ elements and for
$I=\{i_1,\dots,i_k\}$ we define
$I^{\sigma}=\{\sigma(i_1),\dots,\sigma(i_k)\}$. Furthermore, for
each $I\subset\{1,\dots,K\}$ we get the diagram
$$
\begin{array}{rcccl}
&V&\stackrel{\chi}{-\,-\,\to}&V&\\
\pi_I&\downarrow&&\downarrow&\pi_{I^{\sigma}}\\
&F_I&\stackrel{\chi_I}{-\,-\,\to}&F_{I^{\sigma}},&
\end{array}
$$
where $\chi_I$ is a birational map. In particular, $\chi$ induces
birational isomorphisms
$$
\chi_e\colon F_e-\,-\,\to F_{\sigma(e)},
$$
$e=1,\dots,K$. However, all the varieties $F_e$ are birationally
superrigid, so that all the maps $\chi_e$ are biregular
isomorphisms. Thus
$$
\chi=(\chi_1,\dots,\chi_K)\in \mathop{\rm Bir}V
$$
is a biregular isomorphism, too: $\chi\in \mathop{\rm Aut}V$.
Q.E.D. for Corollary 1.


\section{Divisorially canonical varieties}

In this section we prove Theorem 2.

\subsection{Fano hypersurfaces}


Let $F=F_M\subset{\mathbb P}={\mathbb P}^M$ be a smooth Fano
hypersurface. For a point $x\in F$ fix a system of affine
coordinates $z_1,\dots,z_M$ with the origin at $x$ and let
$$
f=q_1+q_2+\dots+q_M
$$
be the equation of the hypersurface $F$, $q_i=q_i(z_*)$ are
homogeneous polynomials of degree $\mathop{\rm deg} q_i=i$. In
order to prove that this variety is divisorially canonical, we
need stronger regularity conditions (that is, the conditions of
general position) for $f$ than those that were used for proving
birational rigidity in [28]. Set
$$
f_i=q_1+\dots+q_i
$$
to be the truncated equation, $i=1,\dots,M$.

($R1.1$) The sequence
$$
q_1,q_2,\dots,q_{M-1}
$$
is regular in ${\cal O}_{x,{\mathbb P}}$, that is, the system of
equations
$$
q_1=q_2=\dots=q_{M-1}=0
$$
defines a one-dimensional subset, a finite set of lines in
${\mathbb P}$, passing through the point $x$. This is the standard
regularity condition, which was used in the previous papers (see
[28]).

($R1.2$) The linear span of any irreducible component of the
closed algebraic set
$$
q_1=q_2=q_3=0
$$
in ${\mathbb C}^M$ is the hyperplane $q_1=0$ (that is, the tangent
hyperplane $T_xF$).

($R1.3$) The closed algebraic set
\begin{equation}\label{b1}
\overline{\{f_1=f_2=0\}\cap F}=\overline{\{q_1=q_2=0\}\cap
F}\subset {\mathbb P}
\end{equation}
(the bar $\bar{}$ means the closure in ${\mathbb P}$) is
irreducible and any section of this set by a hyperplane $P\ni x$
is

\begin{itemize}

\item either also irreducible and reduced,

\item or breaks into two irreducible components
$B_1+B_2$, where $B_i=F\cap S_i$ is the section of $F$ by a plane
$S_i\subset{\mathbb P}$ of codimension 3, and moreover
$\mathop{\rm mult}\nolimits_x B_i=3$,

\item or is non-reduced and is of the form $2B$, where
$B=F\cap S$ is the section of $F$ by a plane $S$ of codimension 3,
and moreover $\mathop{\rm mult}\nolimits_x B=3$.

\end{itemize}


Set ${\cal F}_{\rm reg}\subset {\cal F}$ to be the set of Fano
hypersurfaces, satisfying the conditions ($R1.1-R1.3$) at every
point (in particular, every hypersurface $F\in {\cal F}_{\rm
reg}$ is smooth). It is clear that ${\cal F}_{\rm reg}$ is a
Zariski open subset of the projective space ${\cal F}$.

{\bf Proposition 4.} {\it For $M\geq 6$ the set ${\cal F}_{\rm
reg}$ is non-empty.}

{\bf Proof} is given below in Sec. 2.3. For $M\geq 8$ in the
condition ($R1.3$) we may require that the section of the set
(\ref{b1}) by any hyperplane $P\ni x$ were irreducible and
reduced, a general hypersurface satisfies this stronger
condition. On the other hand, for $M=4,5$ it is easy to show that
for any hypersurface $F\in {\cal F}$ there is a point where the
conditions ($R1.2$) and ($R1.3$) are not satisfied.

Let us prove the claim (i) of Theorem 2, that is, that the
condition ($C$) is satisfied for a regular Fano hypersurface $F\in
{\cal F}_{\rm reg}$.

Let $\Delta\in |nH|$ be an effective divisor, $n\geq 1$, where
$H\in \mathop{\rm Pic}F$ is the class of a hyperplane section,
$K_F=-H$. We have to show that the pair $(F,\frac {1}{n}\Delta)$
has canonical singularities.

Assume the converse. Then for a certain sequence of blow ups
$\varphi\colon F^+\to F$ and an exceptional divisor $E^+\subset
F^+$ the Noether-Fano inequality
\begin{equation}\label{b2}
\nu_{E^+}(\Delta)> n\cdot a(E^+)
\end{equation}
is satisfied. For a fixed $E^+$ the inequality (\ref{b2}) is
linear in $\Delta$, so that without loss of generality we may
assume that $\Delta\subset F$ is a prime divisor, that is, an
irreducible subvariety of codimension 1. From (\ref{b2}) it
follows easily that the centre $Y=\varphi(E^+)$ of the valuation
$E^+$ on $F$ satisfies the inequality
$$
\mathop{\rm mult}\nolimits_Y\Delta> n.
$$
On the other hand, it is well known [25,28], that for any
irreducible curve $C\subset F$ the inequality
$$
\mathop{\rm mult}\nolimits_C\Delta \leq n
$$
holds. Thus $Y=x$ is a point. Let $\varepsilon\colon \widetilde
F\to F$ be its blow up, $E\subset \widetilde F$ the exceptional
divisor $E\cong {\mathbb P}^{M-2}$. By Proposition 3, for some
hyperplane $B\subset E$ the inequality

\begin{equation}
\label{b3} \mathop{\rm mult}\nolimits_x\Delta+\mathop{\rm
mult}\nolimits_B\widetilde\Delta>2n,
\end{equation}
holds, where $\widetilde\Delta\subset\widetilde F$ is the strict
transform of the divisor $\Delta$.

Let ${\mathbb T}=\overline{T_xF}\subset{\mathbb P}$ be the tangent
hyperplane at the point $x$. The divisor $E$ can be naturally
identified with the projectivization
$$
{\mathbb P}(T_x{\mathbb T})={\mathbb P}(T_x F).
$$
There is a unique hyperplane ${\mathbb B}\subset{\mathbb T}$,
$x\in{\mathbb B}$, such that
$$
B={\mathbb P}(T_x{\mathbb B})
$$
with respect to the above-mentioned identification. Let
$\Lambda_{{\mathbb B}}$ be the pencil of hyperplanes in ${\mathbb
P}$, containing ${\mathbb B}$, and $\Lambda_B=\Lambda_{{\mathbb
B}}|_F\subset |H|$ its restriction onto $F$. Consider a general
divisor $R\in\Lambda_B$. It is a hypersurface  of degree $M$ in
${\mathbb P}^{M-1}$, smooth at the point $x$. Let $\widetilde
R\subset\widetilde F$ be the strict transform of the divisor $R$.
Obviously,
$$
\widetilde R\cap E=B.
$$
Set $\Delta_R=\Delta|_R=\Delta\cap R$. It is an effective divisor
on the hypersurface $R$.

{\bf Lemma 3.} {\it The following estimate holds:}
\begin{equation}
\label{b4} \mathop{\rm mult}\nolimits_x \Delta_R>2n.
\end{equation}

{\bf Proof.} We have
$$
(\widetilde \Delta\circ\widetilde R)={\widetilde \Delta}_R+Z,
$$
where $Z$ is an effective divisor on $E$. According to the
elementary rules of the intersection theory [8],
$$
\mathop{\rm mult}\nolimits_x\Delta_R=\mathop{\rm
mult}\nolimits_x\Delta+\mathop{\rm deg} Z,
$$
since $\mathop{\rm mult}\nolimits_x R=1$. However, $Z$ contains
$B$ with multiplicity at least $\mathop{\rm
mult}\nolimits_B\widetilde \Delta$. Therefore, the inequality
(\ref{b3}) implies the estimate (\ref{b4}). Q.E.D. for the lemma.

{\bf Lemma 4.} {\it The divisor
$$
T_R=T_xR\cap R
$$
on the hypersurface $R$ is irreducible and has multiplicity
exactly 2 at the point $x$.}

{\bf Proof.} The irreducibility is obvious (for instance, for
$M\geq 6$ one can apply the Lefschetz theorem). By the condition
($R1.2$) the quadric
$$
\{q_2|_E=0\}
$$
does not contain a hyperplane in $E$ as a component, in
particular, it does not contain the hyperplane $B\subset E$. Thus
the quadratic component of the equation of the divisor $T_R$, that
is, the polynomial
$$
q_2|_B,
$$
is non-zero. Q.E.D. for the lemma.

Let us continue our proof of Theorem 2. By Lemmas 3 and 4 we can
write
$$
\Delta_R=aT_R+\Delta^{\sharp}_R,
$$
where $a\in{\mathbb Z}_+$ and the effective divisor
$\Delta^{\sharp}_R\in |n^{\sharp}H_R|$ on the hypersurface $R$
satisfies the estimate
$$
\mathop{\rm mult}\nolimits_x \Delta^{\sharp}_R>2n^{\sharp}.
$$
Moreover, $\Delta^{\sharp}_R$ does not contain the divisor $T_R$
as a component. Without loss of generality we can assume the
divisor $\Delta^{\sharp}_R$ to be irreducible and reduced.

Now consider the second hypertangent system [28,29,34]
$$
\Lambda^R_2=|s_0f_2+s_1f_1|_R,
$$
where $s_i$ are homogeneous polynomials of degree $i$ in the
linear coordinates $z_{\ast}$. Its base set
$$
S_R=\{q_1|_R=q_2|_R=0\}
$$
is by condition ($R1.3$) of codimension 2 in $R$ and either
irreducible and of multiplicity 6 at the point $x$, or breaks into
two plane sections of $R$, each of multiplicity 3 at the point
$x$. In any case, for a general divisor $D\in \Lambda^R_2$ we get
$\Delta^{\sharp}_R\not\subset \mathop{\rm Supp} D$, so that the
following effective cycle of codimension two on $R$,
$$
\Delta_D=(D\circ \Delta^{\sharp}_R),
$$
is well defined. Since $\mathop{\rm mult}\nolimits_x D=3$ and
$\Lambda^R_2\subset |2H_R|$, the cycle $\Delta_D$ satisfies the
estimate
\begin{equation}
\label{b5} \frac{\mathop{\rm mult}\nolimits_x}{\mathop{\rm
deg}}\Delta_D> \frac{3}{M}.
\end{equation}
We can replace the cycle $\Delta_D$ by its suitable irreducible
component and thus assume it to be an irreducible subvariety of
codimension 2 in $R$. Comparing the estimate (\ref{b5}) with the
description of the set $S_R$ given above, we see that
$$
\Delta_D\not\subset S_R.
$$
This implies that $\Delta_D\not\subset T_R$. Indeed, if this were
not true, we would have got
\begin{equation}
\label{b6} f_1|_{\Delta_D}\equiv q_1|_{\Delta_D}\equiv 0.
\end{equation}
However, $\Delta_D\subset D$, so that for some $s_0\neq 0$,
$s_1\neq 0$ (the divisor $D$ is chosen to be general) we have
$$
(s_0f_2+s_1f_1)|_{\Delta_D}\equiv 0.
$$
By (\ref{b6}) this implies that
$$
f_2|_{\Delta_D}\equiv (q_1+q_2)|_{\Delta_D}\equiv 0
$$
(since $s_0\neq 0$ is just a constant), so that
$$
\Delta_D\subset S_R.
$$
A contradiction.

Thus $\Delta_D\not\subset T_R$. Therefore the effective cycle
$$
\Delta^+=(\Delta_D\circ T_R)
$$
is well defined. It satisfies the estimate
\begin{equation}
\label{b7} \frac{\mathop{\rm mult}\nolimits_x}{\mathop{\rm
deg}}\Delta^+>\frac{6}{M}.
\end{equation}
The effective cycle $\Delta^+$ as a cycle on $F$ is of codimension
4.
Now recall the following fact [28]: if the Fano hypersurface $F$
at the point $x$ satisfies the regularity condition ($R1.1$), then
for any effective cycle $Y$ of pure codimension $l\leq M-2$ the
inequality
$$
\frac{\mathop{\rm mult}\nolimits_x}{\mathop{\rm deg}}Y\leq
\frac{l+2}{M}
$$
holds. Therefore, the inequality (\ref{b7}) for an effective
cycle of codimension 4 is impossible.

The proof of the claim (i) of Theorem 2 is complete.

\subsection{Fano double spaces}
Let $F\stackrel{\sigma}{\to}{\mathbb P}={\mathbb P}^M$ be a Fano
double space branched over a smooth hypersurface $W=W_{2M}\subset
{\mathbb P}$ of degree $2M$, $M\geq 3$. For a point $x\in W$ fix
a system of affine coordinates $z_1,\dots,z_M$ on ${\mathbb P}$
with the origin at $x$ and set
$$
w= q_1+ q_2+\dots + q_{2M}
$$
to be the equation of the hypersurface $W$, $q_i= q_i(z_*)$ are
homogeneous polynomials of degree $\mathop{\rm deg}q_i=i$. The
regularity condition is formulated below in terms of the
polynomials $q_i$. One should consider the three cases $M\geq 5$,
$M=4$ and $M=3$ separately. For convenience of notations assume
that $q_1\equiv z_1$. Set also
$$
\bar q_i=\bar
q_i(z_2,\dots,z_M)=q_i|_{\{z_1=0\}}=q_i(0,z_2,\dots,z_M).
$$
Now let us formulate the regularity condition at the point $x$.

($R2$) Let $M\geq 5$. We require that the rank of the quadratic
form $\bar q_2$ is at least 2.

Assume that $M=3$ or 4. We require that the quadratic form $\bar
q_2$ is non-zero and moreover

(i) either $\mathop{\rm rk}{\bar q}_2\geq 2$ (as above),

(ii) or $\mathop{\rm rk}{\bar q}_2=1$ and the following
additional condition is satisfied. Without loss of generality we
assume in this case that
$$
\bar q_2= z^2_2.
$$
Now for $M=4$ we require that the following cubic polynomial in
the variable $t$,
$$
\bar q_3(0,1,t)= q_3(0,0,1,t)
$$
has three distinct roots.

For $M=3$ we require that at least one of the following two
polynomials in the variable $z_3$,
$$
\bar q_3(0,z_3)\quad \mbox{or}\quad \bar q_4(0,z_3)
$$
(they are of the form $\alpha z^m_3$, $\alpha\in{\mathbb C}$,
$m=3,4$) is non-zero.

Set ${\cal W}_{\rm reg}\subset {\cal W}$ to be the set of Fano
double spaces, satisfying the condition ($R2$) at every point of
the branch divisor. It is obvious that ${\cal W}_{\rm reg}$ is a
Zariski open subset.

{\bf Proposition 5.} {\it The set ${\cal W}_{\rm reg}$ is
non-empty.}

{\bf Proof} is given below in Sec. 2.3.

Now let us prove the claim (ii) of Theorem 2, that is, that a
regular Fano double cover $F$ satisfies the condition ($C$).

Assume the converse: there is an irreducible divisor $D\in |nH|$,
$H=-K_F=\sigma^*H_{\mathbb P}$, where $H_{\mathbb P}$ is the
class of a hyperplane in ${\mathbb P}$, such that the pair
$(F,\frac{1}{n}D)$ is not canonical. Recall the following fact
[26, Proposition 4.3]: for any irreducble curve $C\subset F$ the
inequality
$$
\mathop{\rm mult}\nolimits_CD\leq n
$$
holds. (Formally speaking, this inequality was proved in [26] for
a {\it movable} divisor $D\in \Lambda\subset |nH|$, where
$\Lambda$ is a linear system without fixed components, however the
movability was never used in the proof.) Now arguing as in Sec.
2.1 above, we conclude that the centre of a non-canonical
singularity of the pair $(F,\frac{1}{n}D)$ is a point $x\in F$.
Again we have the inequality
$$
\mathop{\rm mult}\nolimits_xD+\mathop{\rm
mult}\nolimits_B\widetilde D> 2n
$$
for some hyperplane $B\subset E$ in the exceptional divisor
$E\subset \widetilde F$ of the blow up $\widetilde F\to F$ of the
point $x$.


{\bf Case 1.} The point $z=\sigma(x)\not\in W$. In this case let
$P\subset{\mathbb P}$ be the hyperplane such that $P\ni z$ and
$$
{\widetilde F}_P\cap E=B,
$$
where $F_P=\sigma^{-1}(P)$. Such hyperplane is unique. The divisor
$F_P$ is obviously irreducible and satisfies the equality
$$
\mathop{\rm mult}\nolimits_x F_P+\mathop{\rm mult}\nolimits_B
{\widetilde F}_P=2.
$$
Since $F_P\in |H|$, this implies that $F_P\neq D$. Therefore,
$$
Z=(F_P\circ D)
$$
is a well defined effective cycle of codimension two on $F$. We
get
$$
\mathop{\rm deg} Z=2n,\quad \mathop{\rm mult}\nolimits_x Z\geq
\mathop{\rm mult}\nolimits_x D+\mathop{\rm mult}\nolimits_B
\widetilde D> 2n,
$$
which is impossible. Thus the case $z\not\in W$ does not realize.

{\bf Case 2.} The point $z=\sigma(x)\in W$. In this case let
$$
P=T_z W\subset{\mathbb P}
$$
be the tangent hyperplane. Set again $F_P=\sigma^{-1}(P)$.

We get $F_P\in |H|$, $\mathop{\rm mult}\nolimits_x F_P\geq 2$. The
divisor $F_P$ is irreducible, so that if $F_P\neq D$, then for
the cycle $Z=(F_P\circ D)$ we get again
\begin{equation} \label{b7a}
\mathop{\rm deg} Z=2n,\quad \mathop{\rm mult}\nolimits_x Z\geq
2\mathop{\rm mult}\nolimits_x D > 2n,
\end{equation}
since $\mathop{\rm mult}\nolimits_x D\geq \mathop{\rm
mult}\nolimits_B \widetilde D$. Since there is no cycle $Z$,
satisfying the condition (\ref{b7a}), we conclude that
$$
F_P=D.
$$
Up to this moment we have never used the arguments of general
position.


{\bf Lemma 5.}  {\it The pair $(F,F_P)$ is canonical at the point
$x$.}

{\bf Proof} is given below in Sec. 2.3. This contradiction
completes the proof of Theorem 2.

\subsection{The regularity conditions}

Let us prove Proposition 4. We have to show that the closed set
$$
{\cal F}\setminus {\cal F}_{\rm reg}
$$
of Fano hypersurfaces, non-regular at at least one point, is of
positive codimension in ${\cal F}$. Let ${\cal F}(x)$ be the set
of hypersurfaces passing through a point $x\in {\mathbb P}$,
${\cal F}_{\rm reg}(x)\subset {\cal F}(x)$  the set of
hypersurfaces, regular at the point $x$. It suffices to show that
\begin{equation}\label{b8}
\mathop{\rm codim}\nolimits_{{\cal F}(x)}[{\cal F}(x)\setminus
{\cal F}_{\rm reg}(x)]\geq M,
\end{equation}
since ${\cal F}(x)$ is a divisor in ${\cal F}$. For the first
condition ($R1.1$) it has already been proved in [28], so that we
may assume that all hypersurfaces under consideration satisfy
($R1.1$). Fix a point $x\in {\mathbb P}$ and show that violation
of the condition ($R1.2$) or ($R1.3$) imposes at least $M$
conditions on the hypersurface $F\ni x$ (that is, on the
polynomial $f$). As we shall see, the actual codimension of the
set of non-regular hypersurfaces is much higher, and moreover,
the estimate is sharper as $M$ gets higher. Let
$$
\varphi\colon\widetilde F\to F
$$
be the blow up of the point $x$, $E\subset \widetilde F$ the
exceptional divisor, $E\cong {\mathbb P}^{M-2}$. We consider
$z_2,\dots,z_M$ as homogeneous coordinates on $E$, assuming that
$q_1\equiv z_1$. Set
$$
g_i=q_i|_E=q_i|_{\{z_1=0\}}.
$$
Now the condition ($R1.2)$ reads as follows: the linear span of
any irreducible component of the set $\{g_2=g_3=0\}\subset E$ is
$E$.

It is easier to check the inequality (\ref{b8}) as $M$ gets
higher. We consider in full detail the hardest case $M=6$. Here
$E={\mathbb P}^4$. It is easy to see that reducibility of the
quadric
$$
Q=\{g_2=0\}\subset E
$$
imposes on the polynomial $g_2$ (and thus on $q_2$) 6 independent
conditions. Thus for a general hypersurface the quadric $Q$ is
irreducible at every point $x\in F$. Therefore, we get the
following possibilities for $Q$:

\begin{itemize}

\item the quadric $Q$ is smooth (this is true for a point $x\in F$
of general position),

\item the quadric $Q$ is a cone over a smooth quadric in
${\mathbb P}^3$ (this is the case for the points lying on a
divisor $Z_1\subset F$),

\item the quadric $Q$ is a cone with the vertex ${\mathbb P}^1$
over a smooth conic in ${\mathbb P}^2$ (this happens for the
points $x\in Z_2$, where $Z_2\subset F$ is a closed set of
codimension 3).

\end{itemize}

Let us consider each of these three cases separately.

Assume that the quadric $Q$ is smooth. We have
$$
\mathop{\rm Pic}Q={\mathbb Z}H_Q,
$$
where $H_Q$ is the hyperplane section. The condition ($R1.2$) is
violated only in the case when the divisor
$$
\{g_3|_Q=0\}\subset Q
$$
breaks into two components, one of which is a hyperplane section.
It is easy to check that this imposes on the polynomial $g_3$
(and thus on $q_3$) 12 independent conditions, so that for a
general hypersurface $F$ the condition ($R1.2$) is not violated
when the quadric $Q$ is smooth.

Note that for an arbitrary $M\geq 6$ in the case of smooth quadric
$Q$ violation of the condition ($R1.2$) imposes
$$
\frac16(M^3-3M^2-10M+24)
$$
independent conditions on the polynomial $q_3$, so that the extra
codimension of the set of non-regular hypersurfaces increases when
$M$ grows higher.

Now assume that the quadric $Q$ is a cone with the vertex at a
point $p\in E$ or a line $L\subset E$ and some component of the
divisor $\{g_3|_Q=0\}$ is contained in a hyperplane $P\subset E$.
Consider the intersection $P\cap Q$. If this intersection is an
irreducible quadric, then we argue as above in the smooth case. If
$P\cap Q$ is reducible, then
$$
P\cap Q=S_1+S_2
$$
for some planes $S_i\subset Q$. Vanishing of the polynomial $g_3$
on a plane $S\subset Q$ imposes $10-1=9$ independent conditions on
$g_3$ and $q_3$, since the quadric $Q$ contains a one-dimensional
family of planes.

If $M\geq 7$, then the arguments are similar and the estimates
stronger, see the example above. Thus we have proved that the
condition ($R1.2$) holds at every point of a general hypersurface
$F\subset {\mathbb P}$ for $M\geq 6$.

Let us consider the condition ($R1.3$). Let
$$
{\mathbb T}=\overline{\{q_1=0\}}\subset {\mathbb P}
$$
be the projective tangent hyperplane to the hypersurface $F$ at
the point $x$, ${\mathbb T}\cong {\mathbb P}^{M-1}$ and
$$
Q=\{q_2|_{\mathbb T}=0\}\subset {\mathbb T}
$$
a quadratic cone with the vertex at the point $x$. To check the
condition ($R1.3$) we must inspect a number of cases. The
arguments are of the same type, for this reason we consider just
a few main examples. Let us explain the scheme of arguments.

Fix the quadric $Q$ and a hyperplane $P\ni x$, $P\subset {\mathbb
T}$. Set $Q_P=Q\cap P$. We look at the intersection
$$
Q_P\cap F
$$
as a divisor on the quadric $Q_P$. Obviously, the quadric
$Q_P\cap F$ is given by the equation
$$
(q_3+\dots+q_M)|_{Q_P}=0,
$$
where the homogeneous polynomials $q_3,\dots,q_M$ are arbitrary.
Now to prove the proposition, it is necessary to estimate, how
many conditions on the polynomials $q_3,\dots,q_M$ imposes a
violation of ($R1.3$) for each type of the quadric $Q_P$. Again
we restrict ourselves by the hardest case $M=6$.

Assume at first that $Q$ is a cone over a non-singular
three-dimensional quadric. In this case $Q$ is a factorial
variety. Let $P\ni x$ be an arbitrary hyperplane. For the quadric
$Q_P=Q\cap P$ the two cases are possible:

(1)  $Q_P$ is a cone with the vertex $x$ over a non-singular
quadric $S\subset {\mathbb P}^3$;

(2)  $Q_P$ is a cone with the vertex at a line $L\ni x$ over a
non-singular conic $C\subset{\mathbb P}^2$. This possibility
realizes if and only if the hyperplane $P$ is tangent to $Q$
along the line $L$.

In its turn, the case (1) breaks into two subcases.

(1A) The set $Q_P\cap F$ is reducible and at least one of its
components, say $A$, does not contain the point $x$. For the
divisor $A\subset Q_P$ we get
$$
A\sim aH_P,
$$
where $H_P$ is the hyperplane section of the quadric $Q_P$ and
$1\leq a\leq 4$. (A little bit later we explain why $a\leq 4$:
formally $a\leq 5$.) Thus
$$
Q_P\cap F= A+B,
$$
where $B\sim (6-a)H_P$. The affine equation of the hypersurface
$F\cap P\subset P$ takes the form
$$
q^+_2+\dots+q^+_6=0,
$$
where $q^+_i=q_i|_P$: that is why $a\leq 4$. The equation
$q^+_2=0$ defines the quadric $Q_P$. The set of sections of the
quadric $Q_P$ by hypersurfaces $F$ of this form is of dimension
125 whereas for $a\geq 2$
$$
h^0(Q_P,{\cal O}(a))={a+4 \choose 4}- {a+2 \choose 4}.
$$
Now an easy calculation shows that with the hyperplane $P\subset
{\mathbb T}$ fixed violation of the condition ($R1.3$) imposes on
the set of polynomials $(q^+_3,\dots,q^+_6)$ at least 35
independent conditions. Since $P\ni x$ varies in a 4-dimensional
family, we obtain finally that violation of the condition
($R1.3$) of the type ($1A$) at the point $x\in F$ imposes on $f$
at least 31 conditions. Recall that we need just 6 conditions.
Thus for a general hypersurface $F$ the subcase ($1A$) is
impossible.

($1B$) Here the set $Q_P\cap F$ is reducible and each of its
components contains the point $x$. In this case the following
sets are reducible (or everywhere non-reduced): the projectivized
tangent cone
$$
C={\mathbb P}T_x(Q_P\cap F)
$$
and the intersection
$$
C_{\infty}=Q_P\cap F\cap H_{\infty}
$$
with the hyperplane at infinity $H_{\infty}\not\ni x$. Since
$Q_P$ is a cone over the smooth quadric $S\cong Q_P\cap
H_{\infty}$, both curves $C$ and $C_{\infty}$ can be looked at as
curves on the quadric $S\cong {\mathbb P}^1\times {\mathbb P}^1$.
The curve $C$ is given by the polynomial $q_3$, the curve
$C_{\infty}$ by the polynomial $q_6$. Reducibility of the curve
$C$ (which is of bidegree $(3,3)$) imposes $k\geq 3$ independent
conditions on the polynomial $q_3$, and moreover, the value $k=3$
corresponds to the case of general position when $C=R+C^+$, where
$R\subset S$ is a line and $C^+$ an irreducible curve of degree
5. Reducibility of the curve $C_{\infty}$ (which is of bidegree
$(6,6)$) imposes $k_{\infty}\geq 6$ independent conditions on the
polynomial $q_6$, and moreover, the value $k_{\infty}=6$
corresponds to the case of general position when
$C_{\infty}=R_{\infty} +C^+_{\infty}$, where $R_{\infty}\subset
S$ is a line and $C^+_{\infty}$ is an irreducible curve of degree
11. However, by the condition ($R1.1$) the hypersurface $F$ (and
thus $Q_P\cap F$) cannot contain a plane. Therefore a violation
of the condition ($R1.3$) imposes on the polynomial $f$
$$
k+k_{\infty}\geq 3+6+1=10
$$
independent conditions. (Taking into account the polynomials
$q_4$ and $q_5$ gives a considerably higher codimension, but we
do not need that.) Since the hyperplane $P$ varies in a
4-dimensional family ($x\in P\subset {\mathbb T}$), we get
finally at least 6 independent conditions on the polynomial $f$.
Therefore, on a general hypersurface $F$ a violation of the
condition ($R1.3$) of the type (1) is impossible.

Let us consider the case (2). Here the hyperplane $P\subset
{\mathbb T}$ varies in a 3-dimensional family. Assume that
$Q_P\cap F$ is reducible. If $F\not\supset L$, where $L\cong
{\mathbb P}^1$ is the vertex of the cone $Q_P$, then the
arguments completely similar to the case ($1A$) above give the
required estimate (\ref{b8}).

The condition $F\supset L$ means that
$$
q_i|_L\equiv 0,\quad i=3,4,5,6.
$$
which gives 4 additional independent conditions on $f$.
Furthermore, it is easy to check that the case
$$
\mathop{\rm mult}\nolimits_LQ_P\cap F>2
$$
on a general hypersurface $F$ is impossible. Let $S_1,S_2\subset
Q_P$, $S_i\supset L$ be a pair of general planes on $Q_P$. The
restriction $F|_{S_i}$ is a plane curve of degree 6 with the line
$L$ as a component of multiplicity 1. Reducibility of the
residual curves of degree 5 imposes on $f$ at least 6 independent
conditions. As a result, we get at least $7=6+4\, - \, 3$
independent conditions on $f$ at the point $x$. Recall that
$Q_P\cap F$ cannot contain planes by the condition ($R1.1$).

The case when $Q$ is a cone over a non-singular three-dimensional
quadric is completed.

The other cases are inspected in a similar way. We do not give
these arguments because they are of the same type and quite
elementary. Let us consider in detail just one more case which is
opposite to the case (1) above. It occurs for $M=6$ only.

Let $Q$ be a cone with the vertex at a plane $S\cong {\mathbb
P}^2$ over a conic $C\subset {\mathbb P}^2$. Moreover assume that
the hyperplane $P$ contains $S$, so that
$$
Q_P=\Pi_1+\Pi_2\subset P \quad \mbox{or} \quad Q_P=2\Pi
$$
is a pair of 3-planes or a double 3-plane. The set of points $x\in
F$ where this situation is possible is two-dimensional. The
hyperplane $P$ varies in a two-dimensional family: $S\subset
P\subset {\mathbb T}$. The hypersurface $F$ cannot contain planes
because of the condition ($R1.1$).

If $\Pi_1\neq\Pi_2$, then $\Pi_1\cap \Pi_2=S$ and
$$
F\cap Q_P=F\cap\Pi_1+F\cap \Pi_2.
$$
Reducibility of the surface $F\cap \Pi_i$ imposes on $F$ at least
25 conditions so that we may assume that both surfaces $F\cap
\Pi_i$ are irreducible. Obviously,
$$
\mathop{\rm mult}\nolimits_x F\cap \Pi_i\geq 3.
$$
However, by the condition ($R1.2$) we have
$$
\mathop{\rm mult}\nolimits_x F\cap Q_P=\mathop{\rm
mult}\nolimits_xF\cap \Pi_1+\mathop{\rm mult}\nolimits_x F\cap
\Pi_2=6.
$$
Therefore the condition ($R1.3$) holds in this case, too.

If $\Pi_1=\Pi_2=\Pi$, the arguments are similar.

This completes our proof for $M=6$. For $M=7$ the arguments are
simpler and for $M\geq 8$ one can argue like in the case (1)
above.

Q.E.D. for Proposition 4.

{\bf Proof of Proposition 5.} The set of quadratic forms of rank
$\leq 1$ in the variables $z_2,\dots,z_M$ is of codimension
$$
c(M)=\frac{(M-1)(M-2)}{2}.
$$
When $M\geq 5$ we have $c(M)\geq M$, so that ($R2$) holds at
every point $x\in W$ for a sufficiently general polynomial $w$.

Assume that $M=4$. Here $c(4)=3=\mathop{\rm dim}W$, so that for a
general hypersurface $W$ the condition (R2), (i) is violated at a
finite set of points $x\in W$. For a general cubic polynomial
$q_3$ the polynomial $\bar q_3(0,1,t)$ has three distinct roots:
a multiple root gives a condition of codimension 1 for $\bar
q_3$. Therefore ($R2$) holds for $M=4$.

Let $M=3$. Here $c(3)=1$, $\mathop{\rm dim} W=2$, so that there
is a one-dimensional set of points $x\in W$ where the conditon
(i) is violated. It is easy to see that the conditions
$$
\bar q_3(0,z_3)\equiv 0 \quad \mbox{and} \quad \bar
q_4(0,z_3)\equiv 0
$$
are independent.

This completes the proof of Proposition 5.

{\bf Proof of Lemma 5.} Since the hypersurface $W$ is
non-singular, the point $x$ is an isolated singularity of the
hypersurface
$$
W_P=W\cap P,
$$
where $P=T_xW$. The irreducible divisor $F_P$ is the double cover
of the hyperplane $P$ branched at $W_P$. Set $p=\sigma^{-1}(x)\in
F$, $\varphi\colon\widetilde F\to F$ the blow up of the point $p$,
$E\subset \widetilde F$ the exceptional divisor. Assume that the
pair ($F,F_P$) is not canonical. By Proposition 3, for some
hyperplane $B\subset E$ the inequality
\begin{equation}\label{b9}
\mathop{\rm mult}\nolimits_pF_P+\mathop{\rm
mult}\nolimits_B\widetilde F_P\geq 3
\end{equation}
holds. However, $\mathop{\rm mult}\nolimits_p F_P=2$ so that the
projectivized tangent cone
\begin{equation}\label{b10}
{\mathbb P}(T_pF_P)=(\widetilde F_P\circ E)
\end{equation}
must contain the hyperplane $B$. On the other hand, the tangent
cone $T_pF_P$ is given by the equation
$$
y^2=\bar q_2(z_2,\dots,z_M)
$$
(in the notations of Sec. 2.2). For $M\geq 5$ by the condition
($R2$) we get $\mathop{\rm rk}\bar q_2\geq 2$, so that the
quadric (\ref{b10}) is irreducible and cannot contain $B$. This
proves the lemma for $M\geq 5$.

When $M=3,4$ the inequality (\ref{b9}) holds only for those points
$p$ where $\mathop{\rm rk}\bar q_2=1$ and only for the two
hyperplanes
$$
B=\{y\pm z_2=0\}.
$$
Let $M=4$. Local computations show that $\widetilde F_P$ has
exactly three singular points: they lie on the line
$$
\{y=z_2=0\}
$$
and correspond to the roots of the polynomial $\bar q_3(0,1,t)$.
Moreover, these three points are non-degenerate quadratic
singularities. Thus the three-fold $F_P$ has terminal
singularities, so that the pair $(F,F_P)$ is terminal, the more
so canonical, contrary to our assumption.

Let $M=3$. Here $E={\mathbb P}^2$ and
$$
\widetilde F_P\cap E= L_++L_-=\{y=\pm z_2\}
$$
is a pair of lines. The only singularity of the surface
$\widetilde F_P$ over the point $p$ is the point $p^*=L_+\cap L_-$
of intersection of these lines. If $\bar q_3(0,z_3)\neq 0$, then
$p^*$ is a non-degenerate quadratic singularity, so that its blow
up
$$
\varphi^*\colon F^*_P\to \widetilde F_P
$$
gives a non-singular surface $F^*_P$. If $\bar q_3(0,z_3)\equiv
0$, then the exceptional divisor $E^*$ of the blow up $\varphi^*$
is a pair of lines, $E^*=L^*_++L^*_-$, so that the only
singularity of the surface $F^*_P$ over the point $p^*$ is the
point $p^{\sharp}=L^*_+\cap L^*_-$.

It is easy to check that by the condition $\bar q_4(0,z_3)\neq 0$
the point $p^{\sharp}$ is a non-degenerate quadratic singularity.

Thus in any case the surface $F_P$ is canonical.

Q.E.D. for the lemma.


\section{The connectedness principle of Shokurov and Koll\' ar}

In this section we discuss the connectedness principle and
inversion of adjunction and prove Proposition 3.

\subsection{The connectedness principle}
The following fact is a particular case of the general
connectedness principle [15, Theorem 17.4].

{\bf Proposition 6 (Shokurov, Koll\' ar).} {\it Let $x\in X$ be a
smooth germ and
$$
D=\sum_{i\in I}d_iD_i
$$
an effective ${\mathbb Q}$-divisor ($d_i\in {\mathbb Q}_+$ for all
$i\in I$). Let
$$
\varphi\colon \widetilde X\to X
$$
be a resolution of singularities of the pair ($X,D$), in
particular, the support of the divisor $\widetilde D$ is a
divisor with normal crossings. Write down
$$
K_{\widetilde X}= \varphi^*(K_X+D)+\sum_{j\in J}e_jE_j,
$$
where $E_j\subset \widetilde X$ are prime divisors (either
exceptional divisors of the morphism $\varphi$, or components of
the strict transform $\widetilde D$). Then the divisor
$$
\bigcup_{e_j\leq -1}E_j
$$
is connected in a neighborhood of every fiber of the map
$\varphi$.}

{\bf Proof:} see [15, Chapter 17]. Note that the main ingredient
of the proof is the Kawamata-Viehweg vanishing theorem [4,17,43].


\subsection{Inversion of adjunction}

Inversion of adjunction follows from the connectedness principle.
In the general form inversion of adjunction is formulated and
proved in [15, Chapter 17]. We use a particular case of this
general fact.

{\bf Proposition 7 (inversion of adjunction).} {\it Let $x\in X$
be a smooth germ, as above, and $D$ an effective ${\mathbb
Q}$-divisor, $D=\sum\limits_{i\in I}d_iD_i$ with $0\leq d_i\leq
1$, $x\in \mathop{\rm Supp} D$. Assume that the point $x$ is an
isolated centre of a non-canonical singularity of the pair
($X,D$), that is, the pair ($X,D$) is not canonical, but its
restriction onto $X\setminus \{x\}$ is canonical. Let $R\ni x$ be
a smooth divisor, $R\not\subset \mathop{\rm Supp}D$. Then the pair
$$
(R,D_R=D|_R)
$$
is not log canonical at the point $x$.}

{\bf Proof.} For convenience of the reader we give a proof,
following [15, Chapter 17]. Replacing $D$ by
$$
\frac{1}{1+\varepsilon}D
$$
for a small $\varepsilon \in {\mathbb Q}_+$, we may assume that
$d_i<1$ for all $i\in I$. The pair $(X,D)$ remains non-canonical.

Let $\lambda\colon X^+\to X$ be the blow up of the point $x$,
$E=\lambda^{-1}(x)\subset X^+$ the exceptional divisor, $D^+$ and
$R^+$ the strict transforms of the divisors $D$ and $R$,
respectively. Furthermore, let
$$
\mu\colon \widetilde X \to X^+
$$
be a resolution of singularities of the pair ($X^+,D^++ R^+$),
$$
\varphi=\lambda\circ \mu\colon \widetilde X\to X
$$
the composite map. Now write down
\begin{equation}\label{c1}
K_{\widetilde X}=\varphi^*(K_X+D+R)+\sum_{j\in J}e_jE_j-\sum_{i\in
I}d_i\widetilde D_i-\widetilde R,
\end{equation}
where $E_j$, $j\in J$, are all exceptional divisors of the
morphism $\varphi$, $\widetilde D_i$ and $\widetilde R$ are the
strict transforms of the divisors $D_i$ and $R$ on $\widetilde
X$, respectively. Set
$$
b_j=\mathop{\rm ord}\nolimits_{E_j}\varphi^*D,\quad a_j=a(E_j,X),
$$
$j\in J$. In these notations we get for $j\in J$
$$
e_j=a_j-b_j-r_j,
$$
where $r_j=\mathop{\rm ord}_{E_j}\varphi^*R$. Obviously,
$$
\varphi^{-1}(x)=\bigcup_{j\in J^+}E_j
$$
for a subset $J^+\subset J$. If $j\in J^+$, then
$$
r_j=\mathop{\rm ord}\nolimits_{E_j}\mu^*R^++\mathop{\rm
ord}\nolimits_{E_j}\mu^*E\geq 1.
$$
By our assumption the pair ($X,D$) is not canonical, but
canonical outside the point $x$. Therefore, there is an index $l$
among $j\in J^+$ such that $a_l<b_l$. For this index we have
$$
e_l<-1.
$$
By the connectedness principle we may assume that
$$
E_l\cap\widetilde R\neq\emptyset.
$$
Now from (\ref{c1}) by the adjunction formula we get
$$
K_{\widetilde R}=(K_{\widetilde X}+\widetilde R)|_{\widetilde R}=
\varphi^*_R(K_R+D_R)+(\sum_{j\in J}e_jE_j|_{\widetilde
R}-\sum_{i\in I}d_i\widetilde D_i|_{\widetilde R}),
$$
where $\varphi_R=\varphi|_{\widetilde R}\colon \widetilde R\to R$
is the restriction of the sequence of blow ups $\varphi$ onto
$R$. By what has been said in the last bracket there is at least
one prime divisor of the form $E_l|_{\widetilde R}$, where $l\in
J^+$, with the coefficient strictly smaller than $-1$. Q.E.D. for
Proposition 7.

Here is one more corollary from the connectedness principle.

{\bf Proposition 8.} {\it In the notations of the proof of
Proposition 7 assume that the pair ($X,D$) is not log canonical,
but log canonical outside the point $x$. Then the following
alternative takes place:

(1) either $\mathop{\rm mult}_xD>\mathop{\rm dim}X$,

(2) or the set
$$
\mu(\bigcup_{b_j>a_j+1}E_j)\subset E
$$
 is connected.}

{\bf Proof.} By the assumptions the claim follows immediately
from the connectedness principle.


\subsection{Proof of Proposition 3}

Canonicity is stronger than log canonicity. Therefore one can
apply inversion of adjunction (Proposition 7) several times,
subsequently restricting the pair ($X,D$) onto smooth subvarieties
$$
R_1\supset R_2\supset\dots\supset R_k,
$$
where $R_1\subset X$ is a smooth divisor, $R_{i+1}\subset R_i$ is
a smooth divisor, $x\in R_k$ and $R_k\not\subset \mathop{\rm
Supp}D$. All the pairs
$$
(R_i,D|_{R_i})
$$
are not log canonical at the point $x$. Thus Proposition 7 holds
for any smooth germ $R\ni x$ of codimension $k\leq \mathop{\rm
dim}X-1$. In particular, it holds for a general surface $S\ni x$.
(This fact was for the first time used by Corti [1] in order to
obtain an alternative proof of the $4n^2$-inequality, see also
[34, Proposition 1.5].) Thus the pair
$$
(S,D_S=D|_S)
$$
has at the point $x$ an isolated (for a general $S$) non log
canonical singularity. Let us consider the two-dimensional case
more closely. Let
$x\in S$ be a germ of a smooth surface, $C\subset S$ a germ of an
effective (possibly reducible) curve, $x\in C$. Consider a
sequence of blow ups
$$
\varphi_{i,i-1}\colon S_i\to S_{i-1},
$$
$S_0=S$, $i=1,\dots,N$, $\varphi_{i,i-1}$ blows up a point
$x_{i-1}\in S_{i-1}$, $E_i=\varphi^{-1}_{i,i-1}(x_{i-1})\subset
S_i$ is the exceptional line. For $i>j$ set
$$
\varphi_{i,j}=\varphi_{j+1,j}\circ\dots\circ\varphi_{i,i-1}\colon
S_i\to S_j,
$$
$\varphi=\varphi_{N,0}$, $\widetilde S=S_N$. We assume that the
points $x_i$ lie one over another, that is, $x_i\in E_i$, and that
$x_0=x$, so that all the points $x_i$, $i\geq 1$, lie over $x$:
$$
\varphi_{i,0}(x_i)=x\in S.
$$
Let $\Gamma$ be the graph with the vertices $1,\dots,N$ and
oriented edges (arrows)
$$
i\to j,
$$
that connect $i$ and $j$ if and only if $i>j$ and
$$
x_{i-1}\in E^{i-1}_j,
$$
where for a curve $Y\subset S_j$ its strict transform on $S_a$,
$a\geq j$, is denoted by the symbol $Y^a$. Assume that the point
$x$ is the centre of an isolated non log canonical singularity of
the pair
$$
(S,\frac{1}{n}C)
$$
for some $n\geq 1$. This means that for some exceptional divisor
$E\subset \widetilde S$ the inequality
\begin{equation}\label{c2}
\nu_E(C)=\mathop{\rm ord}\nolimits_E\varphi^*C> n(a_E+1),
\end{equation}
holds, where $a_E$ is the discrepancy of $E$. Without loss of
generality we may assume that $E=E_N$ is the last exceptional
divisor. The inequality (\ref{c2}) is called {\it the log
Noether-Fano inequality}.

For $i>j$ let the symbol $p_{ij}$ denote the number of paths in
the graph $\Gamma$ from the vertex $i$ to the vertex $j$. For
$i<j$ set $p_{ij}=0$. Set also $p_{ii}=1$. In terms of the numbers
$p_{ij}$ the log Noether-Fano inequality (\ref{c2}) takes the
traditional form
\begin{equation}\label{c3}
\sum^N_{i=1}p_{Ni}\mu_i>n(\sum^N_{i=1}p_{Ni}+1),
\end{equation}
where
$$
\mu_i=\mathop{\rm mult}\nolimits_{x_{i-1}}C^{i-1}.
$$

{\bf Proposition 9.} {\it Either $\mu_1>2n$ (that is, the first
exceptional divisor $E_1\subset S_1$ already gives a non log
canonical singularity of the pair $(S,(1/n)C)$), or $N\geq 2$ and
the following inequality holds:}
$$
\mu_1+\mu_2>2n.
$$

{\bf Proof.} If $N=1$, then $\mu_1>2n$ by means of log
Noether-Fano inequality. Assume that $\mu_1\leq 2n$, then $N\geq
2$. Obviously, $\mu_1> n$. If $\mu_2\geq n$, then $\mu_1+\mu_2>
2n$, as we claim. So assume that $\mu_2< n$. Then for each
$i\in\{2,\dots,N\}$ we have $\mu_i\leq \mu_2< n$ (since the point
$x_{i-1}$ lies over $x_1$). Therefore from the inequality
(\ref{c3}) we get
$$
p_{N1}(\mu_1-n)+\sum^N_{i=2}p_{Ni}(\mu_2-n)> n.
$$
However,
$$
p_{N1}=\sum_{j\to 1}p_{Nj}\leq \sum^N_{i=2}p_{Ni},
$$
so that the more so
$$
\sum^N_{i=2}p_{Ni}(\mu_1+\mu_2-2n)> n.
$$
Therefore $\mu_1+\mu_2 > 2n$. Q.E.D. for the proposition.


{\bf Proof of Proposition 3.} We prove Proposition 3 in the
following form:

{\it let ($X,D$) be a pair as in Proposition 8, $\lambda\colon
X^+\to X$ the blow up of the point $x$, $E\subset X^+$ the
exceptional divisor, $E\cong {\mathbb P}^e$, $e=\mathop{\rm
dim}X-1$. Then there exists a hyperplane $B\subset E$ such that
\begin{equation}\label{c4}
\mathop{\rm mult}\nolimits_xD+\mathop{\rm mult}\nolimits_BD> 2,
\end{equation}
where $D^+$ is the strict transform of $D$ on $X^+$.}

{\bf Proof.} Consider a general surface $S\ni x$. The pair
$(S,D_S)$ is not log canonical, but log canonical outside the
point $x$. By Proposition 8, either
$$
\mathop{\rm mult}\nolimits_xD_S>2,
$$
but in this case $\mathop{\rm mult}_xD>2$ and the inequality
(\ref{c4}) holds in a trivial way for any hyperplane $B\subset
E$, or the centres of all non log canonical singularities on $S^+$
(that is, the strict transform of $S$ on $X^+$) are contained in
a proper closed connected subset
$$
Z_S\subset E_S=E\cap S^+\cong {\mathbb P}^1.
$$
Obviously, $Z_S$ is a point $y_S\in E_S$. Since the surface $S$
is general, there is a hyperplane $B\subset E$ such that
$$
y_S=B\cap S^+.
$$
By Proposition 9, the inequality
$$
\mathop{\rm mult}\nolimits_x D_S+\mathop{\rm
mult}\nolimits_{y_S}D^+_S>2
$$
holds. This immediately implies the inequality (\ref{c4}) and
Proposition 3.

\newpage

{\small

\section*{References}

\noindent 1. Corti A., Singularities of linear systems and 3-fold
birational geometry, in ``Explicit Birational Geometry of
Threefolds'', London Mathematical Society Lecture Note Series
{\bf 281} (2000), Cambridge University Press, 259-312.
\vspace{0.3cm}

\noindent 2. Corti A., Pukhlikov A. and Reid M., Fano 3-fold
hypersurfaces, in ``Explicit Birational Geometry of Threefolds'',
London Mathematical Society Lecture Note Series {\bf 281} (2000),
Cambridge University Press, 175-258. \vspace{0.3cm}

\noindent 3. Corti A. and Reid M., Foreword to ``Explicit
Birational Geometry of Threefolds'', London Mathematical Society
Lecture Note Series {\bf 281} (2000), Cambridge University Press,
1-20. \vspace{0.3cm}

\noindent 4. Esnault H. and Viehweg E., Lectures on vanishing
theorems, DMV-Seminar. Bd. {\bf 20.} Birkh\" auser, 1992.
\vspace{0.3cm}

\noindent 5. Fano G., Sopra alcune varieta algebriche a tre
dimensioni aventi tutti i generi nulli, Atti Acc. Torino {\bf 43}
(1908), 973-977. \vspace{0.3cm}

\noindent 6. Fano G., Osservazioni sopra alcune varieta non
razionali aventi tutti i generi nulli, Atti Acc. Torino {\bf 50}
(1915), 1067-1072. \vspace{0.3cm}

\noindent 7. Fano G., Nouve ricerche sulle varieta algebriche a
tre dimensioni a curve-sezioni canoniche, Comm. Rent. Ac. Sci.
{\bf 11} (1947), 635-720. \vspace{0.3cm}

\noindent 8. Fulton W., Intersection Theory, Springer-Verlag,
1984. \vspace{0.3cm}

\noindent 9. Graber T., Harris J. and Starr J., Families of
rationally connected varieties. J. Amer. Math. Soc. {\bf 16}
(2002), no. 1, 57-67. \vspace{0.3cm}

\noindent 10. Grinenko M. M., Birational automorphisms of a
three-dimensional double cone, Sbornik: Mathematics {\bf 189}
(1998), no. 7-8, 991-1007. \vspace{0.5cm}

\noindent 11. Grinenko M. M., Birational properties of pencils of
del Pezzo surfaces of degrees 1 and 2, Sbornik: Mathematics {\bf
191} (2000), no. 5-6, 633-653. \vspace{0.5cm}

\noindent 12. Grinenko M. M., Birational properties of pencils of
del Pezzo surfaces of degrees 1 and 2. II, Sbornik: Mathematics
{\bf 194} (2003), no. 5, 669-696. \vspace{0.3cm}

\noindent 13. Iskovskikh V.A., Birational automorphisms of
three-dimensional algebraic varieties, J. Soviet Math. {\bf 13}
(1980), 815-868. \vspace{0.3cm}

\noindent 14. Iskovskikh V.A. and Manin Yu.I., Three-dimensional
quartics and counterexamples to the L\" uroth problem, Math. USSR
Sb. {\bf 86} (1971), no. 1, 140-166. \vspace{0.3cm}

\noindent 15. Koll{\'a}r J., et al., Flips and Abundance for
Algebraic Threefolds, Asterisque 211, 1993. \vspace{0.3cm}

\noindent 16. Koll\' ar J., Rational curves on algebraic
varieties. Springer, 1996. \vspace{0.3cm}

\noindent 17. Kawamata Y., A generalization of Kodaira-Ramanujam's
vanishing theorem, Math. Ann. {\bf 261} (1982),
43-46.\vspace{0.3cm}

\noindent 18. Manin Yu. I., Rational surfaces over perfect fields.
Publ. Math. IHES {\bf 30} (1966), 55-113. \vspace{0.3cm}

\noindent 19. Manin Yu. I. Rational surfaces over perfect fields.
II. Mat. Sb. {\bf 72} (1967), 161-192. \vspace{0.3cm}

\noindent 20. Manin Yu. I., Cubic forms. Algebra, geometry,
arithmetic. Second edition. North-Holland Mathematical Library,
{\bf 4.} North-Holland Publishing Co., Amsterdam, 1986.
\vspace{0.3cm}

\noindent 21. Noether M., {\" U}ber Fl{\" a}chen welche Schaaren
rationaler Curven besitzen, Math. Ann. {\bf 3} (1871), 161-227.
\vspace{0.3cm}

\noindent 22. Pukhlikov A.V., Birational isomorphisms of
four-dimensional quintics, Invent. Math. {\bf 87} (1987),
303-329. \vspace{0.3cm}

\noindent 23. Pukhlikov A.V., Birational automorphisms of a double
space and a double quadric, Math. USSR Izv. {\bf 32} (1989),
233-243. \vspace{0.3cm}

\noindent 24. Pukhlikov A.V., Birational automorphisms of a
three-dimensional quartic with an elementary singularity, Math.
USSR Sb. {\bf 63} (1989), 457-482. \vspace{0.3cm}

\noindent 25. Pukhlikov A.V., A note on the theorem of
V.A.Iskovskikh and Yu.I.Manin on the three-dimensional quartic,
Proc. Steklov Math. Inst. {\bf 208} (1995), 244-254.
\vspace{0.3cm}

\noindent 26. Pukhlikov A.V., Essentials of the method of maximal
singularities, in ``Explicit Birational Geometry of Threefolds'',
London Mathematical Society Lecture Note Series {\bf 281} (2000),
Cambridge University Press, 73-100. \vspace{0.3cm}

\noindent 27. Pukhlikov A.V., Birational automorphisms of
three-dimensional algebraic varieties with a pencil of del Pezzo
surfaces, Izvestiya: Mathematics {\bf 62}:1 (1998), 115-155.
\vspace{0.3cm}

\noindent 28. Pukhlikov A.V., Birational automorphisms of Fano
hypersurfaces, Invent. Math. {\bf 134} (1998), no. 2, 401-426.
\vspace{0.3cm}

\noindent 29. Pukhlikov A.V., Birationally rigid Fano double
hypersurfaces, Sbornik: Mathematics {\bf 191} (2000), No. 6,
101-126. \vspace{0.3cm}

\noindent 30. Pukhlikov A.V., Birationally rigid Fano fibrations,
Izvestiya: Mathematics {\bf 64} (2000), 131-150. \vspace{0.3cm}

\noindent 31. Pukhlikov A.V., Birationally rigid Fano complete
intersections, Crelle J. f\" ur die reine und angew. Math. {\bf
541} (2001), 55-79. \vspace{0.3cm}

\noindent 32. Pukhlikov A.V., Birationally rigid Fano
hypersurfaces with isolated singularities, Sbornik: Mathematics
{\bf 193} (2002), no. 3, 445-471. \vspace{0.3cm}

\noindent 33. Pukhlikov A.V., Birationally rigid Fano
hypersurfaces, Izvestiya: Mathematics {\bf 66} (2002), no. 6,
1243-1269. \vspace{0.3cm}

\noindent 34. Pukhlikov A.V., Birationally rigid iterated Fano
double covers, Izvestiya: Mathematics. {\bf 67} (2003), no. 3,
555-596. \vspace{0.3cm}

\noindent 35. Pukhlikov A.V., Birationally rigid varieties with a
pencil of Fano double covers. I, Sbornik: Mathematics {\bf 195}
(2004), No. 7. \vspace{0.3cm}

\noindent 36. Pukhlikov A.V., Birationally rigid varieties with a
pencil of Fano double covers. II, Sbornik: Mathematics {\bf 195}
(2004), No. 10. \vspace{0.3cm}

\noindent 37. Pukhlikov A.V., Birational geometry of Fano direct
products, arXiv: math.AG/0405011.\vspace{0.3cm}

\noindent 38. Sarkisov V.G., Birational automorphisms of conic
bundles, Math. USSR Izv. {\bf 17} (1981), 177-202. \vspace{0.3cm}

\noindent 39. Sarkisov V.G., On conic bundle structures, Math.
USSR Izv. {\bf 20} (1982), no. 2, 354-390. \vspace{0.3cm}

\noindent 40. Shokurov V.V., 3-fold log flips, Izvestiya:
Mathematics {\bf 40} (1993), 93-202. \vspace{0.3cm}

\noindent 41. Sobolev I. V., On a series of birationally rigid
varieties with a pencil of Fano hypersurfaces. Sbornik:
Mathematics {\bf 192} (2001), no. 9-10, 1543-1551. \vspace{0.3cm}

\noindent 42. Sobolev I. V., Birational automorphisms of a class
of varieties fibered into cubic surfaces. Izvestiya: Mathematics
{\bf 66} (2002), no. 1, 201-222. \vspace{0.3cm}

\noindent 43. Viehweg E., Vanishing theorems, Crelle J. f\" ur die
reine und angew. Math. {\bf 335} (1982), 1-8.

}

\end{document}